\documentclass[12pt,a4paper]{article}

\usepackage{color}
\usepackage{a4wide}
\usepackage{amsfonts,amsmath}
\usepackage{latexsym}
\usepackage{times}
\author{K\'aroly J. B\"or\"oczky and Daniel Hug}
\title{Isotropic measures and stronger forms of the reverse isoperimetric inequality }

\newcommand{\proof}{\noindent{\it Proof: }}
\newcommand{\proofbox}{\mbox{ $\Box$}\\}
\newcommand{\R}{\mathbb{R}}
\newcommand{\N}{\mathbb{N}}

\newcommand{\prob}{\nu}

\newtheorem{lemma}{Lemma}[section]
\newtheorem{theo}[lemma]{Theorem}

\newtheorem{prop}[lemma]{Proposition}

\begin{document}

\maketitle
%\pagebreak

\begin{abstract}
The reverse isoperimetric inequality, due to Keith Ball, states that if $K$ is an $n$-dimensional convex body, then there is an affine image $\tilde{K}$ of $K$ for which
$S(\tilde{K})^n/V(\tilde{K})^{n-1}$  is bounded from above by the corresponding expression for a regular $n$-dimensional simplex, where $S$ and $V$ denote the surface area and volume functional. It was shown by Franck Barthe that the upper bound is attained only if $K$ is a simplex. The discussion of the equality case is based on the equality case in the geometric form of the Brascamp-Lieb inequality. The present paper establishes stability versions of the reverse isoperimetric inequality and of the corresponding
inequality for isotropic measures.
\end{abstract}

{\scriptsize 
\tableofcontents}

\section{Introduction}

The isoperimetric inequality states  that a Euclidean ball has smallest surface area among convex bodies (compact convex sets with non-empty interiors) of given volume in Euclidean space $\R^n$ with scalar product $\langle\cdot\,,\cdot\rangle$ and norm $\|\cdot\|$, and that Euclidean balls are the only minimizers. Let $B^n$ be the Euclidean unit ball centred at the origin. Denoting by $S(K)$ the surface area and by $V(K)$ the volume of a convex body $K$ in $\R^n$, the isoperimetric inequality can be expressed by the inequality
\begin{equation}\label{isop1}
\frac{S(B^n)^n}{V(B^n)^{n-1}}\le \frac{S(K)^n}{V(K)^{n-1}},
\end{equation}
where equality holds if and only if $K$ is a Euclidean ball. Since surface area and volume are continuous functionals (with respect to the Hausdorff metric) and  
the extremal bodies of the inequality \eqref{isop1} are precisely the Euclidean balls, the following question arises naturally. Suppose that a convex body $K$ in $\R^n$ satisfies
$$
\frac{S(K)^n}{V(K)^{n-1}}\le (1+\varepsilon) \frac{S(B^n)^n}{V(B^n)^{n-1}}
$$
for some $\varepsilon\ge 0$. Does it follow that $K$ is $\varepsilon$-close to a Euclidean ball? An answer to this question requires that the distance ${\rm dist}(K)$ of $K$ from a Euclidean ball is measured in a suitable way. For instance, the distance function ${\rm dist}(\cdot)$ should have the same scaling and motion invariance as the isoperimetric problem. 
The problem can also be stated in the following form. Let again $K$ be a convex body in $\R^n$ and assume that ${\rm dist}(K)\ge \varepsilon$ for some $\varepsilon\ge 0$. Does it 
follow that 
$$
\frac{S(K)^n}{V(K)^{n-1}}\ge (1+f(\varepsilon)) \frac{S(B^n)^n}{V(B^n)^{n-1}},
$$
where $f:[0,\infty)\to[0,\infty)$ is a continuous and increasing function with $f(0)=0$? In other words, is it true that 
$$
\frac{S(K)^n}{V(K)^{n-1}}\ge (1+f({\rm dist}(K))) \frac{S(B^n)^n}{V(B^n)^{n-1}} 
$$
with an explicitly given function $f$? Any such inequality  provides a strengthening of the classical isoperimetric inequality and is called a stability result 
related to \eqref{isop1}. 

Although results of this type can be traced back to work of Minkowski and Bonnesen, a systematic exploration  is much more recent. 
Introductory surveys on geometric stability results were given by H.~Groemer \cite{Groemer1990, Groemer1993}, an up-to-date coverage of various aspects (including applications) of the topic is provided throughout R.~Schneider's book \cite{Sch14}. 
More specifically, stability results for the isoperimetric problem (based on 
the Hausdorff distance) have been found, for instance,  by Groemer and Schneider \cite{GroemerSchneider1991}. 
As a recent breakthrough, N.~Fusco, F.~Maggi, A.~Pratelli \cite{FMP08}
obtained an optimal stability version of the isoperimetric
inequality in terms of the volume difference, and
 A.~Figalli, F.~Maggi, A.~Pratelli \cite{FMP09, FMP10}
even extended the result to the
 Brunn-Minkowski inequality.

The ratio ${S(K)^n}/{V(K)^{n-1}}$ is unbounded from above, if $K$ ranges over all convex bodies.
In fact, simple examples show that $K$ can have arbitrarily small volume and still surface area equal to a prescribed positive value. 
In order to avoid this type of situation, it is a well known strategy (see, for instance, F.~Behrend \cite{Behrend1937})  to consider the affine invariant
$$
\mbox{\rm ir}(K):=\inf\left\{\frac{S(\Phi K)^n}{V(\Phi K)^{n-1}}:\Phi\in {\rm GL}(n)\right\}.
$$
The infimum is attained and the unique minimizer can be characterized, as shown by C.~M.~Petty \cite{Petty1961}
(see also A. Giannopoulos, M. Papadimitrakis \cite{GianPapa1999}). In fact, $K$ minimizes the isoperimetric ratio within its affine equivalence class 
if and only if the suitably normalized area measure of $K$ is isotropic (as defined below). As a simple consequence, the regular simplex minimizes the isoperimetric 
ratio within the class of simplices. 
 Since the new functional `$\mbox{\rm ir}$' is affine invariant and upper semi-continuous, it attains its maximum on the space of convex bodies. In the Euclidean plane,
W.~Gustin \cite{Gustin1953} showed that $\mbox{\rm ir}(K)\le \mbox{\rm ir}(T^2)$ with equality if and only if
$K$ is a triangle; here $T^2$ denotes a regular triangle circumscribed about $B^2$. An extension of such a result to higher dimensions turned out to be a formidable problem which resisted its  solution  until K.~M.~Ball \cite{Bal89, Ball91} established
  reverse forms of the isoperimetric inequality. To state one of his main results, note that
$$
V(T^n)=\frac{n^{n/2}(n+1)^{(n+1)/2}}{n!} \qquad\mbox{ and }\qquad S(T^n)=nV(T^n),
$$
where $T^n$ is a regular simplex in $\R^n$ circumscribed about $B^n$.

\bigskip

\noindent{\bf Theorem~A } (K.~M.~Ball) {\it For any convex body $K$ in $\R^n$, there exists some $\Phi\in{\rm GL}(n)$ such that}
$$
\frac{S(\Phi K)^n}{V(\Phi K)^{n-1}}\leq \frac{S(T^n)^n}{V(T^n)^{n-1}}.
$$

\bigskip

It was proved by F.~Barthe \cite{Bar98} that  equality holds in Theorem~A  only if $K$ is a  simplex.

The main objective of this paper is to establish a stability version of the reverse isoperimetric inequality. Following
 \cite{FMP09, FMP10, FMP08}, we define an affine invariant distance
of convex bodies $K$ and $M$ based on the volume difference. For this,  let $\alpha=V(K)^{-{1}/n}$, 
$\beta=V(M)^{-{1}/n}$, and then define 
$$
\delta_{\rm vol}(K,M):=\min\left\{V\left(\Phi (\alpha K)\Delta (x+\beta M)\right):\,
\Phi\in{\rm SL}(n),x\in\R^n\right\}.
$$
We observe that $\delta_{\rm vol}(\cdot,\cdot)$ induces a metric on the affine equivalence classes of convex bodies.

A crucial tool in geometric analysis, and in particular in the proof of the reverse isoperimetric inequality by K.~M.~Ball, is the
 John ellipsoid of a convex body $K$ in $\R^n$. This is the unique ellipsoid of
maximal volume contained in $K$. Obviously, there is an affine image of $K$,
whose John ellipsoid is the Euclidean unit ball $B^n$.  Below (see (\ref{John10}) and (\ref{John20})), we list some properties of the John ellipsoid.
For thorough discussions of the properties of the John ellipsoid, and of convex bodies in general, see K.~M.~Ball \cite{Bal03}, P.~M.~Gruber \cite{Gru07}
or R.~Schneider \cite{Sch14}.

\begin{theo}
\label{vol}
Let $K$ be a convex body in $\R^n$, $n\ge 3$, whose John ellipsoid is a Euclidean ball, and let $\varepsilon\in [0,1)$. 
If $\delta_{\rm vol}(K,T^n)\geq \varepsilon$, then
$$
\frac{S(K)^n}{V(K)^{n-1}}\leq (1-\gamma\varepsilon^4)\frac{S(T^n)^n}{V(T^n)^{n-1}},
$$
where one may choose $\gamma=n^{-250n}$.
\end{theo}

Considering a convex body $K$ which is obtained from $T^n$ by cutting off regular simplices of height $\varepsilon$ at  the vertices of $T^n$ and slabs of width $\varepsilon^{n-1}$ parallel to the facets of $T^n$, one can see that the stability order (the exponent of $\varepsilon$) in Theorem \ref{vol} must be at least $1$. 

In the plane, we obtain a result of optimal stability order.

\begin{theo}
\label{voln=2}
Let $K$ be a convex body in $\R^2$,  whose John ellipsoid is a Euclidean ball, and let $\varepsilon\in [0,1)$. 
If $\delta_{\rm vol}(K,T^n)\geq \varepsilon$, then
$$
\frac{S(K)^2}{V(K)}\leq (1-\gamma\varepsilon)\frac{S(T^2)^2}{V(T^2)},
$$
where one may choose $\gamma=2^{-10}3^{-2}$.
\end{theo}

Theorems \ref{vol} and \ref{voln=2} immediately imply that if $K$ is a convex body in $\R^n$ and $\delta_{\rm vol}(K,T^n)\geq \varepsilon$ for some  $\varepsilon\in [0,1)$, 
then $\text{\rm ir}(K)\le (1-\gamma \varepsilon^4)\ \text{\rm ir}(T^n)$, with  $\gamma$ as in these theorem and with $ \varepsilon^4$ replaced by $\varepsilon$ for $n=2$.

Another affine invariant distance between
convex bodies   is the  Banach-Mazur distance $\delta_{\rm BM}(K,M)$, of convex bodies $K$ and $M$, which is defined by
$$
\delta_{\rm BM}(K,M):=\ln\min\{\lambda\geq 1:\,
K-x\subset \Phi(M-y)\subset \lambda(K-x)\mbox{ for }
\Phi\in{\rm GL}(n),x,y\in\R^n\}.
$$
Again, $\delta_{\rm BM}(\cdot,\cdot)$ induces a metric on the affine equivalence classes of convex bodies. The two metrics are related to each other. It is not 
difficult to see that $\delta_{\rm vol}\le 2e^{n^2}\delta_{\rm BM}$ (see Section \ref{secisoperimetric}). In the reverse direction, we have
$\delta_{\rm BM}\le \gamma \ \delta_{\rm vol}^{\frac{1}{n}}$, where $\gamma$ depends on the dimension $n$ (see \cite[Section 5]{BHenk}), and the exponent $\frac1n$ cannot be replaced by anything larger than $\frac2{n+1}$ as can be seen from the example  of a ball from which a cap is cut off.

\begin{theo}
\label{BM}
Let $K$ be a convex body in $\R^n$ whose John ellipsoid is a Euclidean ball, and let $\varepsilon\in [0,1)$. 
If $\delta_{\rm BM}(K,T^n)\geq \varepsilon$, then
$$
\frac{S(K)^n}{V(K)^{n-1}}\leq (1-\gamma\varepsilon^{\max\{4,n\}})\frac{S(T^n)^n}{V(T^n)^{n-1}},
$$
where one may choose $\gamma=n^{-250n}$.
\end{theo}

Cutting off regular simplices of
edge length $\varepsilon$ at the corners of $T^n$, we see that the  error in Theorem~\ref{BM} can be of order
$\varepsilon^{n-1}$.

In the plane, the aforementioned approach due to W.~Gustin  can be used to establish a stability result of optimal order.

\begin{theo}
\label{2D}
Let $K$ be a convex body in $\R^2$, and let $\varepsilon\in [0,1)$. If $\delta_{\rm BM}(K,T^2)\geq \ \varepsilon$, then 
$$
\text{\rm ir}(K)\le (1-\gamma \varepsilon)\ \text{\rm ir}(T^2),
$$
where we can choose $\gamma=2^{-3}3^{-2}$.
\end{theo}

Since $\delta_{\rm vol}\le 2e^{n^2}\delta_{\rm BM}$, Theorem \ref{2D} implies  for a convex body $K$ in $\R^2$ and $\varepsilon\in [0,1)$ 
that if $\delta_{\rm vol}(K,T^2)\geq \ \varepsilon$, then $\text{\rm ir}(K)\le (1-\gamma \varepsilon)\ \text{\rm ir}(T^2)$, 
where we can choose $\gamma=(2e)^{-4}3^{-2}$. In a different way and with a slightly smaller constant $\gamma$, this is also implied by Theorem \ref{voln=2}.

\medskip

As mentioned before, the proof of the reverse isoperimetric inequality by K.~M.~Ball \cite{Bal89, Ball91} is based on a volume estimate for convex bodies
whose John ellipsoid is the unit ball $B^n$. Let $S^{n-1}$ denote the Euclidean unit sphere.
 According to a classical theorem of F.~John \cite{Joh37} (see also K.~M.~Ball \cite{Bal03}),
$B^n$ is the ellipsoid of maximal volume inside a convex body $K$ if and only if $B^n\subset K$  and there exist
$u_1,\ldots,u_k\in S^{n-1}\cap \partial K$ and $c_1,\ldots,c_k>0$ such that
\begin{align}
\label{John10}
\sum_{i=1}^kc_i u_i\otimes u_i&={\rm Id}_n,\\
\label{John20}
\sum_{i=1}^kc_i u_i&=0,
\end{align}
where ${\rm Id}_n$ denotes the $n\times n$ identity matrix and $\partial K$ is the boundary of $K$.

Following E.~Lutwak, D.~Yang, G.~Zhang \cite{Lutwak0}, let us call
a Borel measure $\mu$ on the unit sphere $S^{n-1}$ isotropic if
$$
\int_{S^{n-1}}u\otimes u\, d\mu(u)={\rm Id}_n.
$$
(All measures in the following are supposed to be Borel measures.) 
In this case, equating traces of both sides we obtain that
\begin{equation}
\label{isotropic-total}
\mu(S^{n-1})=n.
\end{equation}
If, in addition, $\mu$ is centred, that is to say, if
$$\int_{S^{n-1}}u\, d\mu(u)=0,$$
then the origin $0$ is an interior point of the convex hull of the support ${\rm supp}\,\mu$ of $\mu$, and hence
$$
Z(\mu):=\{x\in\R^n:\,\langle x,u\rangle\leq 1 \mbox{ for }u\in{\rm supp}\,\mu\}
$$
is a convex body. 

\medskip

The crucial statement leading to the reverse isoperimetric inequality is the following.

\bigskip

\noindent{\bf Theorem~B } {\it If $\mu$ is a centred, isotropic measure on $S^{n-1}$, then}
\begin{equation}\label{ginq}
V(Z(\mu))\leq V(T^n).
\end{equation}
{\it Equality holds if and only if $Z(\mu)$ is a regular simplex circumscribed about $B^n$.}

\bigskip

For a discrete measure $\mu$, the inequality \eqref{ginq} is due to K.~M.~Ball \cite{Bal89, Ball91}. The  equality case was clarified by F.~Barthe \cite{Bar98}. The case of an arbitrary centred, isotropic measure was treated by F.~Barthe \cite{Bar04} and
E.~Lutwak, D.~Yang, G.~Zhang \cite{Lutwak1}, where \cite{Lutwak1} also characterized the equality case. The measures on $S^{n-1}$ which have an isotropic linear image are characterized by
K.~J.~B\"or\"oczky, E.~Lutwak, D.~Yang and  G.~Zhang \cite{BLYZ14}, building on work of
E.~A.~Carlen,
and D.~Cordero-Erausquin \cite{CCE09},  J.~Bennett, A.~Carbery, M.~Christ and T.~Tao \cite{BCCT08} and B.~Klartag \cite{Kla10}.
We note that isotropic measures on $\R^n$ play a central role in the
KLS conjecture by R.~Kannan, L.~Lov\'asz and   M.~Simonovits
\cite{KLM95}; see, for instance, 
F.~Barthe and  D.~Cordero-Erausquin \cite{BCE13}, O.~Guedon and  E.~Milman \cite{GuM11} and
B.~Klartag \cite{Kla09}.

To state a stability version of Theorem~B,  we define the ``spherical" Hausdorff distance
of compact  sets $X,Y\subset S^{n-1}$ by the formula
$$
\delta_H(X,Y):=\min\left\{\max_{x\in X}\min_{y\in Y}\angle(x,y),\max_{y\in Y}\min_{x\in X}\angle(x,y)\right\},
$$
where $\angle(x,y)$ denotes the geodesic distance of $x,y$ on $S^{n-1}$.
In addition, for $x\in S^{n-1}$, we write $\delta[x]$ to denote the Dirac measure on $S^{n-1}$ supported on $\{x\}$, that is, if $A\subset S^{n-1}$ is a measurable set, then $\delta[x](A)=1$
if $x\in A$ and zero otherwise.
If $S$ is a regular simplex circumscribed about $B^n$ with contact points $v_0,\ldots,v_n\in S^{n-1}$, then we set
$$\mu_S=\sum_{i=0}^n\frac{n}{n+1}\,\delta[v_i].$$
For the total mass of $\mu_S$ we obtain $\mu_S(S^{n-1})=n$ as for $\mu$ in \eqref{isotropic-total}.

\begin{theo}
\label{Zmustab}
Let $\mu$ be a centred, isotropic measure on $S^{n-1}$, $n\ge 3$, and let $\varepsilon\in[0,1)$. If
$$
V(Z(\mu))\ge (1-\varepsilon) V(T^n),
$$
then there exists
a regular simplex $S$ circumscribed about $B^n$ such that
$$\delta_H({\rm supp}\,\mu,{\rm supp}\,\mu_S)\leq \gamma \varepsilon^{1/4},$$
where one may choose $\gamma=n^{70n}$.
\end{theo}

Each of the corresponding $n+1$ spherical balls of radius $n^{65n}\varepsilon^{1/4}$ has $\mu$-measure 
of order $\frac{n}{n+1}+O(\varepsilon^{1/4})$, and hence
the Kantorovich-Monge-Rubinstein (or the Wasserstein distance) of $\mu$ from $\mu_S$ is $O(\varepsilon^{1/4})$ where the  implied constant in $O(\cdot)$ depends only on $n$ 
(see Section \ref{secisotropic}). 

\medskip

Again we obtain a result of optimal order for $n=2$. 

\begin{theo}
\label{planariso}
Let $\mu$ be a centred, isotropic measure on $S^1$. If 
$$V(Z(\mu))\geq (1-\varepsilon)V(T^2)$$ 
for $\varepsilon\in[0,1)$, then
there exists a regular triangle $S$ circumscribed about $B^2$ such that
$$
\delta_H({\rm supp}\,\mu,{\rm supp}\,\mu_S)\leq 32\varepsilon.
$$
\end{theo}

\bigskip

We note that the proof of Theorem~B  is based on the rank one case of the geometric Brascamp-Lieb inequality. While we do not actually use the Brascamp-Lieb inequality, an essential tool in  our approach is the proof provided
by  F.~Barthe \cite{Bar97}, which is based on mass transportation. Therefore, it is instructive
to review the argument from \cite{Bar97}, which is done in Section~\ref{secbrascamp-lieb}. At the end of that section,  we outline the arguments leading to
Theorem~\ref{vol}, Theorem \ref{BM} and
Theorem~\ref{Zmustab}  and roughly describe the structure of the paper.

\section{A brief review of the Brascamp-Lieb inequality}
\label{secbrascamp-lieb}

The rank one geometric Brascamp-Lieb inequality, identified by  K. Ball \cite{Bal89}
as an essential case
of the rank one Brascamp-Lieb inequality,
due to H.~J.~Brascamp, E.~H.~Lieb \cite{BrL76},
reads as follows.
If $u_1,\ldots,u_k\in S^{n-1}$ are distinct unit vectors and $c_1,\ldots,c_k>0$ satisfy
$$
\sum_{i=1}^kc_i u_i\otimes u_i={\rm Id}_n,
$$
and $f_1,\ldots,f_k$ are non-negative measurable functions on $\R$, then
\begin{equation}
\label{BL}
\int_{\R^n}\prod_{i=1}^kf_i(\langle x,u_i\rangle)^{c_i}\,dx\leq
\prod_{i=1}^k\left(\int_{\R}f_i\right)^{c_i}.
\end{equation}
According to  F.~Barthe \cite{Bar98}, if equality holds in (\ref{BL}) and none of the functions $f_i$ is identically zero or a scaled version of a  Gaussian, then $k=n$ and  $u_1,\ldots,u_n $ is  an  orthonormal basis of $\R^n$. Conversely, equality holds in (\ref{BL})  if each $f_i$ is a scaled version of the same centered Gaussian,
or if $k=n$ and $u_1,\ldots,u_n$ form an  orthonormal basis.

A thorough discussion of the rank one Brascamp-Lieb inequality can be found in 
E.~Carlen, D.~Cordero-Erausquin \cite{CCE09}. The higher rank case,
due to E.~H.~Lieb \cite{Lie90}, is reproved and further explored by F.~Barthe \cite{Bar98} (including a discussion of the equality case), and is again carefully analysed by
J.~Bennett, T.~Carbery, M.~Christ, T.~Tao \cite{BCCT08}. In particular, see F.~Barthe, D.~Cordero-Erausquin, M.~Ledoux, 
B.~Maurey \cite{BCLM11} for an enlightening review of the relevant literature and an approach via Markov semigroups in a quite general framework. 

 F.~Barthe \cite{Bar97, Bar98}  provides a concise proof of (\ref{BL}) 
based on mass transportation (see also K.~M.~Ball \cite{Bal03}). We sketch the main ideas of this  approach, since this will be the starting point for subsequent refinements.

We assume that
each of the functions $f_i$ is a positive and continuous probability density. Let $g(t)=e^{-\pi t^2}$ be the Gaussian density.
For $i=1,\ldots,k$, we consider the transportation map
$T_i:\R\to\R$ satisfying
$$
\int_{-\infty}^tf_i(s)\,ds=\int_{-\infty}^{T_i(t)} g(s)\,ds.
$$
It is easy to see that $T_i$ is bijective, differentiable and
\begin{equation}\label{masstrans}
f_i(t)=g(T_i(t))\cdot T'_i(t),\qquad t\in\R.
\end{equation}
To these transportation maps, we associate the transformation $\Theta:\R^n\to\R^n$ with
$$
%\label{Thetadef}
\Theta(x):=\sum_{i=1}^kc_iT_i(\langle u_i,x\rangle )\,u_i,\qquad x\in\R^n,
$$
which satisfies
$$
d\Theta(x)=\sum_{i=1}^kc_iT'_i(\langle u_i,x\rangle )\,u_i\otimes u_i.
$$
In this case, $d\Theta$ is positive definite and $\Theta:\R^n\to\R^n$ is
 injective (see \cite{Bar97}). We will need the following two estimates due to K.~M.~Ball \cite{Bal89}.
\begin{description}
\item{(i)} For any $t_1,\ldots,t_k>0$, we have
$$
\det \left(\sum_{i=1}^kt_ic_i u_i\otimes u_i\right)\geq \prod_{i=1}^k t_i^{c_i};
$$
(see also Lemma \ref{BallBarthe} below).
\item{(ii)}
If $z=\sum_{i=1}^kc_i\theta_i u_i$ for
$\theta_1,\ldots,\theta_k\in\R$, then
\begin{equation}
\label{RBLfunc}
\|z\|^2\le \sum_{i=1}^kc_i\theta_i^2.
\end{equation}
\end{description}
Therefore, using first \eqref{masstrans}, and then (i) and  (ii), we obtain
\begin{align*}
\int_{\R^n}\prod_{i=1}^kf_i(\langle u_i,x\rangle)^{c_i}\,dx&=
\int_{\R^n}\left(\prod_{i=1}^kg(T_i(\langle u_i,x\rangle))^{c_i}\right)
\left(\prod_{i=1}^kT'_i(\langle u_i,x\rangle)^{c_i}\right)\,dx\\
&\leq \int_{\R^n}\left(\prod_{i=1}^ke^{-\pi c_iT_i(\langle u_i,x\rangle)^2}\right)
 \det\left(\sum_{i=1}^kc_iT'_i(\langle u_i,x\rangle )\,u_i\otimes u_i\right)\,dx\\
&\leq  \int_{\R^n}e^{-\pi \|\Theta(x)\|^2}\det\left( d\Theta(x)\right)\,dx\\
&\leq  \int_{\R^n}e^{-\pi \|y\|^2}\,dy=1.
\end{align*}

We observe that (i) shows that the optimal constant in the geometric Brascamp-Lieb inequality is $1$.
The stability version  of (i) (with $v_i=\sqrt{c_i} u_i$), Lemma~\ref{Ball-Barthe-stab}, is an essential tool in proving
a stability version of the Brascamp-Lieb inequality leading to Theorem~\ref{Zmustab}.

Let us briefly discuss how K.~M.~Ball \cite{Bal89} used the Brascamp-Lieb inequality to prove
 the discrete version of Theorem~B, since this type of argument is hidden in the proof of
 Proposition~\ref{volZstab} which is crucial for our approach.
 First,  $\R^n$ is embedded into $\R^{n+1}$, and we write $e_{n+1}$
to denote the unit vector in $\R^{n+1}$ orthogonal to $\R^n$.
Let ${\rm supp}\,\mu=\{u_1,\ldots,u_k\}$, let $c_i=\mu(\{u_i\})$, and let
$$
\tilde{u}_i:=-\sqrt{\frac{n}{n+1}} \,u_i+\sqrt{\frac{1}{n+1}}\,e_{n+1}\in S^n 
 \quad\text{  for } i=1,\ldots,k.
$$
The conditions that $\mu$ is isotropic and its centroid is the origin ensure that
$$
\sum_{i=1}^k\tilde{c}_i  \tilde{u}_i\otimes\tilde{u}_i={\rm Id}_{n+1},
\mbox{ \  where $\tilde{c}_i:=\frac{n+1}n\,c_i$ for  $i=1,\ldots,k$.}
$$
Now the Brascamp-Lieb inequality is applied to the system
$\tilde{u}_1,\ldots,\tilde{u}_k,\tilde{c}_1,\ldots,\tilde{c}_k$, where each $f_i$ is the exponential density, that is, $f_i(t)=e^{-t}$ if $t\geq 0$, and $f_i(t)=0$ otherwise.
For the open convex cone
$C=\{y\in\R^{n+1}:\,\langle y,\tilde{u}_i\rangle> 0,\;i=1,\ldots,k\}$, the formulas
(\ref{witness1}) and (\ref{witness2}) in  Section~\ref{secCirc} yield
$$
\int_{\R^{n+1}}\prod_{i=1}^kf_i(\langle y,\tilde{u}_i\rangle)^{\tilde{c}_i}\,dy=
\int_C\exp\left(-\sum_{i=1}^k \tilde{c}_i\langle y,\tilde{u}_i\rangle
  \right)\,dy=V(Z(\mu)) V(T^n)^{-1}.
$$
Since the Brascamp-Lieb inequality implies that this expression is at most $1$, we conclude Theorem~B.

Equality in Theorem~B leads to equality in the Brascamp-Lieb inequality, and hence $k=n+1$ and
$\tilde{u}_1,\ldots,\tilde{u}_{n+1}$ form an orthonormal basis in $\R^{n+1}$. In turn,
$u_1,\ldots,u_{n+1}$ are the vertices of  a regular simplex.

To obtain a stability version of Theorem~B, we need a stability version of the Brascamp-Lieb inequality
in the special case we use. For example, we strengthen (i) in  Section~\ref{secOptcons}, and
 estimate derivatives of the corresponding transportation map
in Section~\ref{secTrans}. The estimates in Section~\ref{secTrans} are very specific for our particular choice o the functions $f_i$, and
no method is known to the authors that could lead to a stability version of the Brascamp-Lieb inequality (\ref{BL}) in general.

The overall structure of the paper is as follows. Sections~\ref{secJohn}, \ref{secOptcons} and
\ref{secalmostreg} provide various  important analytic and geometric estimates concerning John's theorem, related to discrete, isotropic measures and geometric stability results 
for polytopes close to a regular simplex. 
In Section~\ref{secTrans}, we provide auxiliary estimates for the transportation map between the exponential and the Gaussian distribution. After these preparations, 
we establish in Section~\ref{secCirc} the core statement, Proposition~\ref{volZstab}, on which Theorem~\ref{vol}, Theorem \ref{BM} and
Theorem~\ref{Zmustab} are based. Then,  
Section~\ref{secisoperimetric} contains the proofs of   Theorem~\ref{vol} and Theorem~\ref{BM}. In Section \ref{sec2DProof}, we derive  
Theoren \ref{2D}, whose proof is independent of the remaining results. Then, we extend Proposition~\ref{volZstab} to general centred, isotropic measures 
 in Section~\ref{secisotropic}, which proves Theorem~\ref{Zmustab}. Finally, we establish Theorem \ref{planariso} in Section~\ref{2dcase} and Theorem \ref{voln=2} in Section \ref{sec12}.

\section{Some consequences of John's condition}
\label{secJohn}

According to the classical theorem of F.~John \cite{Joh37}, if
$B^n$ is the ellipsoid of maximal volume inside a convex body $K$, then there exist
$u_1,\ldots,u_k\in S^{n-1}\cap \partial K$ and $c_1,\ldots,c_k>0$ such that \eqref{John10} and \eqref{John20} 
are satisfied. 
%\begin{eqnarray}
%\label{John1}
%\sum_{i=1}^kc_i u_i\otimes u_i&=&{\rm Id}_n,\\
%\label{John2}
%\sum_{i=1}^kc_i u_i&=&0.
%\end{eqnarray}
Equating the traces on the two sides of 
\eqref{John10} we obtain 
\begin{equation}
\label{Johntrace}
\sum_{i=1}^k c_i=n.
\end{equation}
In addition, we may assume that 
$$
%\label{Johnbound}
n+1\leq k\leq n(n+3)/2,
$$
where the lower bound on $k$ follows from (\ref{John10}) and  (\ref{John20}) 
and the upper bound on $k$ is implied by the proof of John's theorem \cite{Joh37} (see also  P.~M.~Gruber, F.~E.~Schuster \cite{GrS05}).
We note that (\ref{John10}) is equivalent to
$$
%\label{John1u0}
\sum_{i=1}^k c_i\langle x,u_i\rangle^2=\|x\|^2
\mbox{ \ \ for all $x\in\R^n$}.
$$
Applying this to $x=u_i$ shows that
\begin{equation}
\label{cismall}
 c_i\leq 1
\mbox{ \ \ for $i=1,\ldots,k$}.
\end{equation}

In this section, we discuss properties that only use (\ref{John10}).
This can be written as
\begin{equation}
\label{Johnv}
\sum_{i=1}^k v_i \otimes v_i={\rm Id}_n \mbox{ \ for $v_i:=\sqrt{c_i}\,u_i$}.
\end{equation}
We note that (\ref{Johnv}) is equivalent to
\begin{equation}
\label{John1v0}
\sum_{i=1}^k \langle x,v_i\rangle^2=\|x\|^2
\mbox{ \ \ for all $x\in\R^n$}.
\end{equation}

Given $v_1,\ldots,v_k\in \R^n$ and $\lambda_1,\ldots,\lambda_k>0$,
we consider the $n\times k$ matrix
$$
U:=[\sqrt{\lambda_1}\,v_1,\ldots,\sqrt{\lambda_k}\,v_k].
$$
According to the Cauchy-Binet formula, we have
\begin{equation}
\label{Cauchy-Binet}
\det\left(\sum_{i=1}^k \lambda_i v_i\otimes v_i\right)=\det \left(UU^\top\right)=
\sum_{1\leq i_1<\ldots<i_n\leq k}
\det[\sqrt{\lambda_{i_1}}\,v_{i_1},\ldots,\sqrt{\lambda_{i_n}}\,v_{i_n}]^2.
\end{equation}
It has been pointed out by K.~M.~Ball that the special case $\lambda_1=\ldots=\lambda_k=1$ yields the following estimate.

\begin{lemma}
\label{ciuibig}
If $v_1,\ldots,v_k\in \R^n$
satisfy $\sum_{i=1}^k v_i \otimes v_i={\rm Id}_n$, then
there exist $1\leq i_1<\ldots<i_n\leq k$ such that
$$
\det[v_{i_1},\ldots,v_{i_n}]^2\geq {k \choose n}^{-1}.
$$
\end{lemma}

For non-zero vectors $v$ and $w$, we write
$\angle(v,w)$ to denote their angle, that is, the geodesic distance of the unit vectors $\|v\|^{-1}v$ and $\|w\|^{-1}w$ on the unit sphere.

\begin{lemma}
\label{almostort}
Let  $v_1,\ldots,v_k\in \R^n\setminus \{0\}$
satisfy $\sum_{i=1}^k v_i \otimes v_i={\rm Id}_n$, and let
$0<\eta<1/(3\sqrt{k})$. Assume for any $i\in\{1,\ldots,k\}$ that  $\|v_i\|\leq \eta$ or
there is some $j\in\{1,\ldots,n\}$ with $\angle(v_i,v_j)\leq \eta$.
Then there exists an orthonormal basis $w_1,\ldots,w_n$
such that $\angle(v_i,w_i)<3\sqrt{k}\,\eta$ for $i=1,\ldots,n$.
\end{lemma}
\proof
For $i=1,\ldots,n$, let $u_i=v_i/\|v_i\|$.
We partition the index set $\{1,\ldots,k\}$ into sets
${\cal V}_0,{\cal V}_1,\ldots,{\cal V}_n$
such that  $i\in{\cal V}_i$ for $i=1,\ldots,n$, and in such a way that if $j\in {\cal V}_0$, then $\|v_j\|\leq\eta$,
and if $j\in{\cal V}_i$ for some $i\in\{1,\ldots,n\}$, then $\angle(v_i,v_j)\leq \eta$. 
Observe that ${\cal V}_0$ is possibly empty. 
For $i=1,\ldots,n$, (\ref{John1v0}) yields
$$
1=\|u_i\|^2\geq
\sum_{j\in{\cal V}_i} \langle u_i,v_j\rangle^2\geq
\sum_{j\in{\cal V}_i} \|v_j\|^2\cos^2\eta,
$$
and hence
\begin{equation}
\label{Vi}
\sum_{j\in{\cal V}_i} \|v_j\|^2\leq (\cos\eta)^{-2}.
\end{equation}

For $i=1,\ldots,n$, let $\tilde{w}_i\in S^{n-1}$
 be orthogonal to $v_j$, $j\in\{1,\ldots,n\}\setminus\{i\}$,
and satisfy $\langle \tilde{w}_i,v_i\rangle\geq 0$.
In addition, let $\alpha_i\leq\pi/2$ be the minimal angle
of $\tilde{w}_i$ and any $v_j$ with $j\in{\cal V}_i$, and hence
\begin{equation}
\label{alphai+}
\angle(\tilde{w}_i,v_i)\leq \alpha_i+\eta.
\end{equation}

To bound $\alpha_i$ from above, for $i=1,\ldots,n$, we observe that
$|\langle \tilde{w}_i,v_j\rangle|\leq \eta$  if $j\in {\cal V}_0$. 
Moreover, if $j\in\mathcal{V}_i$, then $\langle \tilde{w}_i,v_j\rangle \le \cos\alpha_i$, and 
if  $j\in\mathcal{V}_l$ for some $l\in\{1,\ldots,n\}\setminus\{i\}$, then $\angle (\tilde{w}_i,v_j)\ge (\pi/2)-\eta$ and therefore 
$\langle \tilde{w}_i,v_j\rangle \le \sin\eta$. 
Using these facts and (\ref{Vi}), we deduce
\begin{align*}
\sum_{j\in{\cal V}_0} \langle \tilde{w}_i,v_j\rangle^2
&\leq (k-n)\eta^2\le \frac{(k-n)\sin^2\eta}{\cos^2\eta},\\
\sum_{j\in{\cal V}_l} \langle \tilde{w}_i,v_j\rangle^2
&\leq\sin^2\eta\sum_{j\in{\cal V}_l} \|v_j\|^2\leq
\frac{\sin^2\eta}{\cos^2\eta},\qquad \mbox{for $l\in\{1,\ldots,n\}\setminus\{i\}$},\\
\sum_{j\in{\cal V}_i} \langle \tilde{w}_i,v_j\rangle^2&\leq
\cos^2\alpha_i\sum_{j\in{\cal V}_i} \|v_j\|^2\leq
\frac{\cos^2\alpha_i}{\cos^2\eta},
\end{align*}
where the sum for ${\cal V}_0$ is set to be zero if ${\cal V}_0$ is empty.
We conclude by (\ref{John1v0}) that
$$
1=\| \tilde{w}_i\|^2\leq  \frac{(k-n)\sin^2\eta}{\cos^2\eta}+
 \frac{(n-1)\sin^2\eta}{\cos^2\eta}+\frac{\cos^2\alpha_i}{\cos^2\eta},
$$
and hence
$$
\sin^2\alpha_i=1-\cos^2\alpha_i\leq 1-\cos^2\eta+(k-1)\sin^2\eta=k\sin^2\eta.
$$
Moreover, for $\eta<1/(3\sqrt{k})$, we have
$$
\frac{\sin(2\sqrt{k}\ \eta)}{\sqrt{k}\sin(\eta)}\ge \frac{\sin(2\sqrt{k}\ \eta)}{\sqrt{k} \ \eta}\ge 2\frac{\sin(2/3)}{2/3}\ge 1.
$$
Therefore,  (\ref{alphai+}) and $ \eta<1/(3\sqrt{k})$ yield
$$
%\label{alphai}
\angle(\tilde{w}_i,v_i)\leq \alpha_i+\eta\leq 2\sqrt{k}\,\eta+\eta<3\sqrt{k}\,\eta, \qquad i=1,\ldots,n.
$$
In particular, this shows that $v_1,\ldots,v_n$ are linearly independent.

We define $w_1=u_1$, and for $i=2,\ldots,n$ we let $w_i$ be
the unit vector in ${\rm lin}\,\{v_1,\ldots,v_i\}$
which is orthogonal to $v_1,\ldots,v_{i-1}$ and satisfies
$\langle w_i,v_i\rangle>0$. Writing
$L_i$ for the orthogonal complement
of ${\rm lin}\,\{v_1,\ldots,v_{i-1}\}$, we have $\tilde{w}_i\in L_i$. 
Since $w_i$ is parallel to the orthogonal projection
of $v_i$ to $L_i$, we conclude that 
$\angle(w_i,v_i)\leq\angle(\tilde{w}_i,v_i)< 3\sqrt{k}\,\eta$.
\proofbox

\section{Analytic stability estimates} 
\label{secOptcons}

To calculate the optimal constant in the
Brascamp-Lieb inequality (\ref{BL}),  the following statement has been
proved by K.~M.~Ball \cite{Bal89}, see  F.~Barthe \cite[Proposition~9]{Bar98} for a simple argument.

\begin{lemma}[K.~M.~Ball]
\label{BallBarthe}
If  $v_1,\ldots,v_k\in \R^n$ satisfy $\sum_{i=1}^kv_i\otimes v_i={\rm Id}_n$ and if \ $t_1,\ldots,t_k>0$, then
$$
\det\left( \sum_{i=1}^kt_iv_i\otimes v_i\right)\geq \prod_{i=1}^k t_i^{\langle v_i,v_i\rangle}.
$$
\end{lemma}
{\bf Remark }  E.~Lutwak, D.~Yang, G.~Zhang \cite{Lutwak0} generalized Lemma~\ref{BallBarthe} for any 
 isotropic measure $\mu$ on $S^{n-1}$ and for any positive continuous function $t$ on ${\rm supp}\,\mu$
in the form
$$
\det\left(\int_{S^{n-1}} t(u)\ u\otimes u\, d\mu(u)\right)\geq
\exp\left(\int_{S^{n-1}}\log t(u)\, d\mu(u)\right),
$$
where equality holds if and only if the quantity $t(v_1)\cdots  t(v_n)$ is constant for linearly independent
$v_1,\ldots,v_n\in {\rm supp}\,\mu$.
Actually  Lemma~\ref{BallBarthe} is the case when
${\rm supp}\,\mu=\{u_1,\ldots,u_k\}$, and $v_i=\sqrt{c_i}\,u_i$ for $c_i=\mu(\{u_i\})$. 
We do not need this generalized version in the present paper.

\bigskip

In Lemma \ref{Ball-Barthe-stab}, we prove a (stronger) stability version of Lemma~\ref{BallBarthe}
by replacing the arithmetic-geometric mean inequality
with the following stability version in the argument of \cite{Bar98}.

\begin{lemma}
\label{geoarithmstab}
If $\prob$ is a probability measure and
$f$ is a  measurable function which is bounded from above and from below by positive constants, then
$$
\frac{\int f\, d\prob}{\exp\left\{\int\ln f \, d\prob\right\}}\ge 1+
\frac{1}{2}\ 
\int\left(\frac{\sqrt{f}}{\sqrt{\int f\, d\prob}}-1\right)^2\, d\prob.
$$
\end{lemma}
\proof  We note that for $a,b\ge 0$, we have
\begin{equation}
\label{Kober}
\frac{a+b}2-\sqrt{a}\ \sqrt{b}=\frac12\left(\sqrt{a}-\sqrt{b}\right)^2.
\end{equation}
Here we choose $b=1$ and
$$
a=\frac{f}{\int f\, d\prob}.
$$
Integrating \eqref{Kober} with this choice of $a,b$ against $\prob$, we get
$$
1-\frac{\int \sqrt{f}\, d\prob}{\sqrt{\int f\, d\prob}}=\frac12
\int\left(\frac{\sqrt{f}}{\sqrt{\int f\, d\prob}}-1\right)^2\, d\prob.
$$
Since $1-x\ge 1-\sqrt{x}$ for $x\in[0,1]$, we obtain
$$
1-\frac{\left(\int \sqrt{f}\, d\prob\right)^2}{\int f\, d\prob}
\ge\frac12 \int\left(\frac{\sqrt{f}}{\sqrt{\int f\, d\prob}}-1\right)^2\, d\prob.
$$
Jensen's inequality yields
$$
\left(\int \sqrt{f}\, d\prob\right)^2\ge \exp\left\{\int\ln f \, d\prob\right\},
$$
and hence we conclude Lemma~\ref{geoarithmstab} by observing that $ (d/c)-1\ge 1-(c/d)$ for any $c,d>0$.
\proofbox

\begin{lemma}
\label{Ball-Barthe-stab}
Let $k\geq n+1$, $t_1,\ldots,t_k>0$, and let $v_1,\ldots,v_k\in \R^n$ satisfy $\sum_{i=1}^kv_i\otimes v_i={\rm Id}_n$. Then
$$
\det\left( \sum_{i=1}^kt_iv_i\otimes v_i\right)\geq \theta^*\ \prod_{i=1}^k t_i^{\langle  v_i,v_i\rangle}
$$
where
\begin{align*}
\theta^*&=  1+\frac12\sum_{1\leq i_1<\ldots<i_n\leq k}
\det[v_{i_1},\ldots,v_{i_n}]^2
\left(\frac{\sqrt{t_{i_1}\cdots  t_{i_n}}}{t_0}-1\right)^2,\\
t_0&=\sqrt{\sum_{1\leq i_1<\ldots<i_n\leq k}
t_{i_1}\cdots t_{i_n}\det[v_{i_1},\ldots,v_{i_n}]^2}.
\end{align*}
\end{lemma}
\proof In this argument, $I$ always denotes some subset of $\{1,\ldots,k\}$ of cardinality $n$. For $I=\{i_1,\ldots,i_n\}$, we define
$$
d_I:=  \det[v_{i_1},\ldots,v_{i_n}]^2\qquad \text{and}\qquad 
t_I := t_{i_1}\cdots t_{i_n}.
$$
From $\sum_{i=1}^kv_i\otimes v_i={\rm Id}_n$ and 
 (\ref{Cauchy-Binet}) we obtain
$$
\sum_I d_I= 1\qquad \text{and}\qquad 
\det \left(\sum_{i=1}^kt_iv_i\otimes v_i\right) =\sum_I t_Id_I,
$$
where the summations extend over all sets $I\subset\{1,\ldots,k\}$ of cardinality $n$. 
It follows that the discrete measure $\mu$ on the
$n$ element subsets of $\{1,\ldots,k\}$ defined by $\mu(\{I\})=d_I$
is a probability measure. According to Lemma~\ref{geoarithmstab},
writing $t_0=\sqrt{\sum_It_Id_I}$, we deduce that
\begin{equation}
\label{agstab}
\det\left( \sum_{i=1}^kt_iv_i\otimes v_i\right) =\sum_I t_Id_I\geq
\left(1+\frac12\sum_Id_I\left( \frac{\sqrt{t_I}}{t_0}-1\right)^2 \right)
\prod_It_I^{d_I}.
\end{equation}
The factor $t_i$ is used in $\prod_It_I^{d_I}$ exactly
$\sum_{I,\,i\in I}d_I$ times. Moreover, (\ref{Cauchy-Binet}) applied
to the vectors $v_1,\ldots,v_{i-1},v_{i+1},\ldots,v_k$ implies
\begin{align*}
\sum_{I,\,i\in I}d_I&=\sum_Id_I-\sum_{I,\,i\not\in I}d_I=
1-\det\left(\sum_{j\neq i}v_j\otimes v_j\right)\\
&=1-\det\left({\rm Id}_n-v_i\otimes v_i\right)=\langle v_i,v_i\rangle.
\end{align*}
Substituting this into (\ref{agstab}) yields the lemma.
\proofbox

To estimate from below (in the proof of Lemma 7.2) the factor $\theta^*$ in Lemma~\ref{Ball-Barthe-stab},
we use the following observation.

\begin{lemma}
\label{xab}
If $a,b,x>0$, then
$$
(xa-1)^2+(xb-1)^2\geq \frac{(a^2-b^2)^2}{2(a^2+b^2)^2}
$$
\end{lemma}
\proof Differentiating $f(x)=(xa-1)^2+(xb-1)^2$ for fixed $a,b$ with
respect to $x$ shows that $f$ attains its minimum at $x=\frac{a+b}{a^2+b^2}$.
Thus
$$
(xa-1)^2+(xb-1)^2\geq \frac{(a-b)^2}{a^2+b^2}=
\frac{(a^2-b^2)^2}{(a^2+b^2)(a+b)^2}\geq
\frac{(a^2-b^2)^2}{2(a^2+b^2)^2}.\mbox{ \ \proofbox}
$$

\section{Polytopes close to  a regular simplex}
\label{secalmostreg}

We prove two quantitative statements about the approximation of a polytope by a simplex. First, we provide 
a lemma which will allow us to put a given orthonormal basis into a more convenient position by a small rotation. 

\begin{lemma}
\label{basis}
Let $e\in S^{n-1}$, and let $\tau\in (0,1/(2n))$. If $w_1,\ldots,w_n$ is an orthonormal basis
of $\R^n$ such that
$$\frac1{\sqrt{n}}-\tau<\langle e,w_i\rangle<\frac1{\sqrt{n}}+\tau\qquad\text{
for $i=1,\ldots,n,$}
$$
then
there exists an orthonormal basis $\tilde{w}_1,\ldots,\tilde{w}_n$
such that $\langle e,\tilde{w}_i\rangle=\frac1{\sqrt{n}}$
and $\angle(w_i,\tilde{w}_i)<n\tau$ for $i=1,\ldots,n$.
\end{lemma}
\proof For $i=1,\ldots,n$, let
$$
\langle e,w_i\rangle=\frac1{\sqrt{n}}+\alpha_i,
\mbox{ \ and hence $|\alpha_i|< \tau$}.
$$
It follows that
$$
1=\|e\|^2=\sum_{i=1}^n\left(\frac1{\sqrt{n}}+\alpha_i\right)^2
<1+\frac2{\sqrt{n}}\left(\sum_{i=1}^n\alpha_i\right)+n\tau^2,
$$
which in turn yields that
$$
\left\langle e,\sum_{i=1}^n w_i\right\rangle=\sqrt{n}+\sum_{i=1}^n\alpha_i>
\sqrt{n}-\frac{n\sqrt{n}}2\,\tau^2>\left\|\sum_{i=1}^n w_i\right\|\cos (n\tau),
$$
since $\cos(n\tau)\le 1-\frac{1}{2}n\tau^2$ for $\tau\in (0,1/(2n))$ and $n\ge 2$. 
In particular, $\angle(e,\sum_{i=1}^n w_i)<n\tau$. We
define $\tilde{w}_i=\Phi (w_i)$ for $i=1,\ldots,n$, where $\Phi$ is the
orthogonal transformation, which rotates $\sum_{i=1}^n w_i$ into
$\sqrt{n}\,e$ via their acute angle in
the two-dimensional linear subspace $L$ containing them, and fixing all vectors in $L^\bot$. Then $\langle e,\tilde w_i\rangle =\langle 
\Phi^{-1} (e),w_i\rangle =\sqrt{n}^{-1}\langle \sum_{j=1}^nw_j,w_i\rangle=1/\sqrt{n}$ for $i=1,\ldots,n$. 
\proofbox

For convex bodies containing the origin in their interiors, we introduce a very specific distance  from
regular simplices whose centroid is the origin. If $K$ is a convex body with
$0\in {\rm int}\,K$, then we define
$$
d(K):=\ln\min\{\lambda\geq 1:\,
sT^n\subset \Phi K\subset \lambda sT^n\mbox{ \ \ for $s>0$ and }
\Phi\in{\rm O}(n)\}.
$$
Clearly, $d(K)=0$ if and only if $K$ is a regular simplex with centroid at the origin.

\begin{lemma}
\label{closesimplex}
Let $Z$ be a polytope, and let $S$ be a regular simplex circumscribed about $B^n$.
Assume that the facets of $Z$ and $S$ touch $B^n$ at $u_1,\ldots,u_k$
and $w_1,\ldots,w_{n+1}$, respectively. Fix $\eta\in (0,1/(9n))$. If 
for any $i\in\{1,\ldots,k\}$ there exists some $j\in\{1,\ldots,n+1\}$ such that $\angle (u_i,w_j)\le \eta$,
then
$$
(1-3n \eta) S\subset  Z\subset (1+3n \eta) S.
$$
In particular, $d(Z)<9n\eta$.
\end{lemma}
\proof The lemma follows from the following
statement: If $\angle (u_1,w_1)\le \eta$ then the tangent plane to $B^n$ at $u_1$
contains $-\lambda w_2$, where
\begin{equation}
\label{lambdaest}
(1-3n \eta)n \le\lambda \le (1+3n \eta)n.
\end{equation}
In order to prove this assertion, we observe that $\lambda^{-1}=\cos\angle (-w_2,u_1)$. Moreover, we write
$\angle (-w_2,u_1)=\alpha+\beta$,
where  $\alpha=\angle (-w_2,w_1)$ with $\cos\alpha=1/n$ and $\tan\alpha<n$,
and $|\beta|\le\eta$. Since
$$
|\cos \beta-1 -\tan\alpha\sin\beta|\le \frac{1}{2}\eta^2 +n\eta=(n+1)\eta
$$
and
$$
|\cos \beta-\tan\alpha\sin\beta|\ge 1-\frac{1}{2}\eta^2-n\eta=1-(n+\eta/2)\eta\ge \frac{1}{2},
$$
we obtain
$$
\left|1-\frac{\lambda}n\right|=
\left|1-\left(\frac{\cos(\alpha+\beta)}{\cos\alpha}\right)^{-1}\right|=
\left|1-\left(\cos\beta-\tan\alpha \sin\beta \right)^{-1}\right|\le 2(n+1)\eta,
$$
which in turn yields (\ref{lambdaest}). 

To conclude the proof, we first observe that the vertices of $S$ are
$-nw_1,\ldots,-nw_{n+1}$. To verify the left inclusion, let $H^-(u):=\{x\in\R^n:\langle x,u\rangle\le1\}$ for $u\in S^{n-1}$. 
We have shown that $-\lambda w_i\in H^-(u_1)$ for $i\in\{2,\ldots,k\}$, and trivially this also 
holds for $i=1$. Hence, \eqref{lambdaest} yields that $(1-3n\eta)(-nw_i)\subset H^-(u_1)$, and therefore $(1-3n\eta) S\subset H^-(u_1)$. Repeating this argument for 
$u_2,\ldots,u_k$, we obtain  $(1-rn\eta) S\subset Z$.

As to the right inclusion, let $tv\in Z$, where $v\in S^{n-1}$ and $t>0$. We can assume that $v$ is in the positive hull of $-w_2,\ldots,-w_{n+1}$. 
Then there is some $i\in \{1,\ldots,k\}$ such that $\angle(u_i,w_1)\le\eta$. 
 By \eqref{lambdaest}, for $j=2,\ldots,n+1$ there are $t_j\in (0,(1+3n\eta)n)$ such that 
$\langle  u_i,-t_jw_j\rangle =1$. There are $\alpha_r\ge 0$ such that $tv=\alpha_2(-w_2)+\ldots+\alpha_{n+1}(-w_{n+1})$, and therefore 
\begin{equation}\label{convexcomb}
\langle u_i,tv\rangle=\sum_{j=2}^{n+1}\langle u_i, {t_j}^{-1}\alpha_j(-t_jw_j)\rangle =\sum_{j=2}^{n+1}\frac{\alpha_j}{t_j}.
\end{equation}
In particular, this shows that  $\langle u_i,v\rangle >0$. Since $tv\in Z$, it is sufficient to prove 
that $tv\in (1+3n\eta)S$ in the case where $\langle u_i,tv\rangle =1$.
But then \eqref{convexcomb} implies that
$$
tv=\sum_{j=2}^{n+1}\frac{\alpha_j}{t_j}(-t_jw_j)\in\text{\rm conv}\{-t_2w_2,\ldots,-t_{n+1}w_{n+1}\}\subset(1+3n\eta)S,
$$
and hence $Z\subset(1+3n\eta)S$. 
\proofbox

\begin{lemma}
\label{closefarsimplex}
Let $Z$ be a polytope, and let $S$ be a regular simplex circumscribed about $B^n$. Fix $\gamma=9\cdot 2^{n+2}n^{2n+2}$ and $\eta\in (0,\gamma^{-1})$. 
Assume that the facets of $Z$ and $S$ touch $B^n$ at $u_1,\ldots,u_k$
and $w_1,\ldots,w_{n+1}$, respectively.  If $\angle (u_i,w_i)\le\eta$ for 
$i=1,\ldots,n+1$  and $\angle (u_k,w_i)\geq \gamma\eta$ for  $i=1,\ldots,n+1$,
then
$$
V(Z)\leq \left(1-\frac{\min_{i=1,\ldots,n+1}\angle (u_k,w_i) }{2^{n+2}n^{2n}}\right)V(S).
$$
\end{lemma}
\proof
Let $H^+:=\{x\in\R^n:\langle x,u_k\rangle\geq 1\}$, and let $F_i$  be the facet of $S$ touching $B^n$ at $w_i$. We may assume that $\angle (u_k,w_1)\leq\angle (u_k,w_i)$ for $i\geq 2$, and hence $\langle u_k,w_1\rangle>0$.

 First, we estimate $V(S\cap H^+)$.
 Let
$z$ be the closest point of $H^+\cap F_1$ to $w_1$. In particular, we have $\|z-w_1\|\leq 1$, while $F_1$ contains the $(n-1)$-ball of radius
$\sqrt{\frac{n+1}{n-1}}>1+\frac1n$ centered at $w_1$. Thus $F_1\cap H^+$ contains a regular $(n-1)$-simplex of height $\frac1n$, and in turn a congruent copy of $\frac1{2n^2}\, F_1$.
 In addition, the  distance of $w_1$ from any $F_i$, $i\geq 2$, is $1+\frac1n$, thus  the  distance of $z$ from  $F_i$ is at least
 $$\frac{1/n}{\|z-w_1\|+(1/n)}\left(1+\frac1n\right)>\frac{h}{2n^2},$$
 where $h=n+1$ is the height of $S$. We deduce that
$H^+\cap S$ contains a point whose distance from $F_1$ is at least
 $\frac{h}{2n^2}\sin\angle (u_k,w_1)$, and hence
$$
V(S\cap H^+)\geq \left(\frac1{2n^2}\right)^{n-1}\  \frac{\angle (u_k,w_1)}{4n^2}\  V(S)=
\frac{\angle (u_k,w_1) }{2^{n+1}n^{2n}}\ V(S).
$$

Let $Z_0$ be the simplex whose facets touch $B^n$ at $u_1,\ldots,u_{n+1}$. Hence
$$(1-3n \eta) S\subset  Z_0\subset (1+3n \eta) S$$
by Lemma~\ref{closesimplex}. It follows that
\begin{align*}
V(Z_0\cap H^+)&\geq V(S\cap H^+)-\left(V(S)-V((1-3n \eta) S)\right)\\
&\ge
\frac{\angle (u_k,w_1) }{2^{n+1}n^{2n}}\  V(S)-3n^2 \eta\  V(S).
\end{align*}
Since $(1+3n\eta)^n\le 1+6n^2\eta$, we have
\begin{align*}
V(Z)&\leq V(Z_0)-V(Z_0\cap H^+)\\
&\le V((1+3n \eta) S)-
\left(\frac{\angle (u_k,w_1) }{2^{n+1}n^{2n}}-3n^2 \eta\right) V(S)\\
&\le \left(1+9n^2 \eta-\frac{\angle (u_k,w_1) }{2^{n+1}n^{2n}}\right) V(S)\\
&\le \left(1-\frac{\angle (u_k,w_1) }{2^{n+2}n^{2n}}\right) V(S),
\end{align*}
which completes the proof. \proofbox

\section{The transportation map}
\label{secTrans}

The argument of F. Barthe \cite{Bar97} uses the transportation
map $\varphi:(0,\infty)\to\R$ between the exponential and the
standard Gaussian density, and hence
\begin{equation}
\label{phidef0}
1-e^{-t}=\int_0^te^{-s}\,ds=\frac1{\sqrt{\pi}}\int_{-\infty}^{\varphi(t)}e^{-s^2}\,ds.
\end{equation}
Clearly, $\varphi$ is strictly increasing and $\varphi(\ln 2)=0$.

\begin{lemma}
\label{phi}
If $t\geq 4$, then $\sqrt{2}<\varphi(t)<\sqrt{t}$, $\frac1{3 \sqrt{t}}<\varphi'(t)<1$ and
 $\varphi''(t)<-\frac1{12 t^{3/2}}$.
\end{lemma}
\proof The definition (\ref{phidef0}) of $\varphi$ can be written in the form
\begin{equation}
\label{phidef}
e^{-t}=\frac1{\sqrt{\pi}}\int_{\varphi(t)}^{\infty}e^{-s^2}\,ds.
\end{equation}
According to the Gordon-Mill inequality (or Mill's ratio, see R.~D.~Gordon \cite{Gor41},  L.~D\"umbgen \cite[(2)]{Dum10}, 
or by a straightforward direct argument), if $z>0$, then
\begin{equation}
\label{GordonMill}
\frac{e^{-z^2}}{2\sqrt{\pi}z}\cdot \frac{2z^2}{2z^2+1}<
\frac1{\sqrt{\pi}}\int_z^{\infty}e^{-s^2}\,ds<\frac{e^{-z^2}}{2\sqrt{\pi}z}.
\end{equation}
We deduce from the left-hand side of (\ref{GordonMill}) that
$$e^{-4}<\frac1{\sqrt{\pi}}\int_{\sqrt{2}}^{\infty}e^{-s^2}\,ds,$$
which in turn implies $\varphi(4)>\sqrt{2}$ by (\ref{phidef}). 
  From (\ref{phidef}) and
the right-hand side of (\ref{GordonMill}), we deduce that $\varphi(t)<\sqrt{t}$ for $t>4$.

We turn to the estimation of derivatives. Differentiating (\ref{phidef}), we get 
\begin{equation}
\label{phidiff}
e^{-t}=\frac{e^{-\varphi(t)^2}\varphi'(t)}{\sqrt{\pi}},\qquad t>0.
\end{equation}
In particular, this shows that $\varphi'(t)>0$ for $t>0$. 
Equation \eqref{phidiff} combined with the right-hand side of  (\ref{GordonMill}) leads to
\begin{equation}
\label{phiphi}
2\varphi(t)\varphi'(t)<1 \mbox{ \ for $t>\ln 2$.}
\end{equation}
Taking the logarithm of (\ref{phidiff}), we deduce  the formula
\begin{equation}
\label{phidifflog}
-t=-\log\sqrt{\pi}-\varphi(t)^2+\log \varphi'(t),
\end{equation}
and differentiating this implies
\begin{equation}
%\label{phi2diff0}
\varphi''(t)=\varphi'(t)(2\varphi(t)\varphi'(t)-1).
\end{equation}
Therefore $\varphi''(t)<0$ follows on the one hand from $\varphi'(t)>0$, and on the other hand
from $\varphi(t)\leq 0$ if $t\leq \ln 2$, and from 
(\ref{phiphi}) if $t>\ln 2$. Thus $\varphi'(t)<\varphi'(\ln 2)=\sqrt{\pi}/2<1$ by (\ref{phidiff})
for $t>\ln 2$.

We also estimate $\varphi''$ in terms of $\varphi$. To this end, we use an improved
version of the right-hand side of the Gordon-Mill inequality  (\ref{GordonMill}) 
(see L. D\"umbgen \cite[(2)]{Dum10}, or by a simple direct argument); namely
$$
\frac1{\sqrt{\pi}}\int_z^{\infty}e^{-s^2}\,ds<
\frac{e^{-z^2}}{2\sqrt{\pi}z}\cdot \frac{2z^2+2}{2z^2+3},\qquad z>0.
$$
We deduce from this and the left-hand side of  (\ref{GordonMill}) that  if $z\geq \sqrt{2}$, then
$$
\frac{e^{-z^2}}{3\sqrt{\pi}z}<\frac1{\sqrt{\pi}}\int_z^{\infty}e^{-s^2}\,ds<
\frac{e^{-z^2}}{2\sqrt{\pi}z}\left(1-\frac1{4z^2}\right).
$$
If $t>4$, then $\varphi(t)>\sqrt{2}$, thus
\begin{equation}
\label{GordonMillgen}
 \frac1{3\varphi(t)}<\varphi'(t)=\sqrt{\pi}e^{\varphi(t)^2-t}< \frac1{2\varphi(t)}
\left(1-\frac1{4\varphi(t)^2}\right).
\end{equation}
In particular, $\varphi'(t) >\frac1{3 \sqrt{t}}$, and combining (\ref{phiphi}) and (\ref{GordonMillgen}) yields
\begin{equation}
%\label{phi2diff}
\varphi''(t)=\varphi'(t)(2\varphi(t)\varphi'(t)-1)<
-\frac{\varphi'(t)}{4\varphi(t)^2}<\frac{-1}{12\varphi(t)^3}
\mbox{ \ for $t>4$, \ }
\end{equation}
which completes the argument. \proofbox

\section{Circumscribed polytopes}
\label{secCirc}

F.~Barthe  \cite{Bar97} proves the Brascamp-Lieb inequality
for functions in one variable in full generality.
This section is based on K.~M.~Ball's \cite{Bal03} interpretation 
 of F.~Barthe's argument in the special case needed for the geometric
application. Since our stability argument uses in an
essential way that the Brascamp-Lieb inequality is
required only for the exponential density function, we do
not separate the statement of the  Brascamp-Lieb inequality.

Proposition \ref{volZstab} is the main ingredient for the proofs of Theorem~\ref{vol}, Theorem \ref{BM} and
Theorem~\ref{Zmustab}. We recall that if $K$ is a convex body with $0\in{\rm int}\ K$, then
$d(K)$ is the minimal $\lambda$ such that there exists a regular simplex $S$ whose centroid is the origin and
$S\subset K\subset e^\lambda S$. 

In the following, we use the abbreviation $N:=n(n+3)/2$. In this section, we consider the case $n\ge 3$, although (with slightly different constants) 
the proof extends also to the case $n=2$. In the plane, however, we can argue in a different way to obtain results of optimal order. For this reason 
we defer the two-dimensional case to Section \ref{2dcase}.

\begin{prop}
\label{volZstab}
Let $\mu$ be a discrete, centred, isotropic measure on $S^{n-1}$. Let $n\ge 3$. Assume that the cardinality of ${\rm supp}\,\mu$ is at most $N+1$, and let 
 $\tau\in(0,n^{-240n})$. If
$$
V(Z(\mu))>(1-\tau)V(T^n),
$$
then there exists a regular simplex $S$ circumscribed about $B^n$  such that
$$
\delta_H({\rm supp}\,\mu,{\rm supp}\,\mu_S)<n^{60n}\tau^{1/4}
\mbox{ \ and \ }d(Z(\mu))<n^{60n}\tau^{1/4}.
$$
\end{prop}

Before we prove Proposition~\ref{volZstab}, we first set up the corresponding notions
following K.~M.~Ball \cite{Bal89}, \cite{Ball91}, and then prove the preparatory statement Lemma~\ref{case12}.

 Let ${\rm supp}\,\mu=\{u_1,\ldots,u_k\}$,
and  let $c_i=\mu(\{u_i\})$. Then 
$\sum_{i=1}^kc_i u_i\otimes u_i={\rm Id}_n$, $\sum_{i=1}^kc_i u_i=0$ and $k\leq N+1$.

We now embed $\R^n$ into $\R^n\times\{0\}=\R^{n+1}$
and write $e_{n+1}$ for the unit vector in $\R^{n+1}$ orthogonal to $\R^n$.
We define
$$
\tilde{u}_i:=-\sqrt{\frac{n}{n+1}}\,u_i+\sqrt{\frac{1}{n+1}}\,e_{n+1}\in S^n\qquad
  \text{and} \qquad  \tilde{c}_i:=\frac{n+1}n\, c_i\qquad \text{for } i=1,\ldots,k,
$$
and hence
\begin{align}
%\label{John1n+1}
\sum_{i=1}^k\tilde{c}_i\ \tilde{u}_i\otimes\tilde{u}_i&={\rm Id}_{n+1},\nonumber\\
\label{John2n+1}
\sum_{i=1}^k\tilde{c}_i\tilde{u}_i&=\sqrt{n+1}\,e_{n+1}, \\
\label{John3n+1}
\sum_{i=1}^k\tilde{c}_i &=n+1.
\end{align}
We observe that if $Z(\mu)$ is a regular simplex circumscribed about $B^n$, then
$k=n+1$ and $ \tilde{u}_1,\ldots, \tilde{u}_{n+1}$ are  an
orthonormal basis of $\R^{n+1}$.

Next we consider the open cone
\begin{align}
\label{Cdef1}
C&:=\{y\in\R^{n+1}:\,\langle y,\tilde{u}_i\rangle> 0,\;i=1,\ldots,k\}\\
\label{Cdef2}
&=\{x+re_{n+1}\in\R^{n+1}:\,x\in\R^n,\;r>0,\;
\langle x,u_i\rangle < r/\sqrt{n},\;i=1,\ldots,k\}
\end{align}
and the map $\Theta:\,C\to\R^{n+1}$ defined by
$$
\Theta(y):=\sum_{i=1}^k  \tilde{c}_i\ \varphi(\langle y,\tilde{u}_i\rangle)\ \tilde{u}_i,
$$
where $\langle y,\tilde{u}_i\rangle>0$ by  (\ref{Cdef1}).
In particular, the differential of $\Theta$ is
$$
%\label{diffT}
d\Theta(y)=\sum_{i=1}^k \tilde{c}_i \ \varphi'(\langle y,\tilde{u}_i\rangle)
\  \tilde{u}_i\otimes\tilde{u}_i .
$$
We observe that $d\Theta$ is positive definite since $\varphi'$ is positive and
$$
\langle z, d\Theta(y)z\rangle=\sum_{i=1}^k \tilde{c}_i\ \varphi'(\langle y,\tilde{u}_i\rangle) \
\langle z,\tilde{u}_i\rangle^2.
$$
It follows that $\Theta$ is injective. 

From (\ref{Cdef2}) we conclude that the section
$\{y\in C:\,\langle y,e_{n+1}\rangle=r\}$ of $C$ for $r>0$ is a translate of
${\rm int}((r/\sqrt{n}) Z(\mu))$. Therefore
\begin{align}
\label{witness1}
\int_Ce^{-\langle y,\sqrt{n+1}\,e_{n+1}\rangle}\,dy&=\int_0^\infty\int_{\frac{r}{\sqrt{n}}\,Z}e^{-\sqrt{n+1}\ r}\,dx\,dr\\
\nonumber
&= V(Z(\mu))\int_0^\infty\left(\frac{r}{\sqrt{n}}\right)^ne^{-\sqrt{n+1}\ r}dr\\
\nonumber
&= V(Z(\mu)) V(T^n)^{-1}.
\end{align}
By first applying (\ref{John2n+1}), then (\ref{phidifflog}),
and finally (\ref{John3n+1}), we deduce  that
\begin{align}
\label{witness2}
\int_Ce^{-\langle y,\sqrt{n+1}\,e_{n+1}\rangle}\,dy&=\int_C\exp\left(-\sum_{i=1}^k \tilde{c}_i\langle y,\tilde{u}_i\rangle
  \right)\,dy\\
\nonumber
&=\int_C\exp\left(\sum_{i=1}^k \tilde{c}_i
(-\log\sqrt{\pi}-\varphi(\langle y,\tilde{u}_i\rangle)^2+\log \varphi'(\langle y,\tilde{u}_i\rangle)
  \right)\,dy\\
\label{Cvol}
&=\pi^{-\frac{n+1}2}\int_C
\exp\left(-\sum_{i=1}^k \tilde{c}_i\varphi(\langle y,\tilde{u}_i\rangle)^2  \right)
\prod_{i=1}^k\varphi'(\langle y,\tilde{u}_i\rangle)^{\tilde{c}_i}\,dy.
\end{align}
For each fixed $y\in C$, we estimate the product of the two
terms in (\ref{Cvol}) after the integral sign.

To estimate the first term  in (\ref{Cvol}), we apply (\ref{RBLfunc})  with $\theta_i=\varphi(\langle y,\tilde{u}_i\rangle)$, and hence the definition
of $\Theta$ yields
\begin{equation}
\label{firstterm}
\exp\left(-\sum_{i=1}^k \tilde{c}_i\varphi(\langle y,\tilde{u}_i\rangle)^2  \right)\leq
 \exp\left(-\|\Theta(y)\|^2\right).
\end{equation}
To estimate the second term, we apply Lemma~\ref{Ball-Barthe-stab}
with  $v_i=\sqrt{\tilde{c}_i}\ \tilde{u}_i$
and $t_i=\varphi'(\langle y,\tilde{u}_i\rangle)$, and
write $\theta(y)$ and $t_0(y)$ to denote
the corresponding $\theta^*\geq 1$ and $t_0$. In particular,
\begin{align}
\label{thetay}
\theta(y)&=  1+\frac12\sum_{1\leq i_1<\ldots<i_{n+1}\leq k}
\tilde{c}_{i_1}\cdots\tilde{c}_{i_{n+1}}
\det[\tilde{u}_{i_1},\ldots,\tilde{u}_{i_{n+1}}]^2\nonumber\\
&\qquad\qquad\times
\left(\frac{\sqrt{\varphi'(\langle y,\tilde{u}_{i_1}\rangle)\cdots
\varphi'(\langle y,\tilde{u}_{i_{n+1}}\rangle)}}{t_0(y)}-1\right)^2,
\end{align}
and  Lemma~\ref{Ball-Barthe-stab} yields
\begin{equation}
\label{secondterm}
\prod_{i=1}^k\varphi'(\langle y,\tilde{u}_i\rangle)^{\tilde{c}_i} \leq
\theta(y)^{-1} \det \left(d\Theta(y)\right).
\end{equation}
We conclude that
\begin{align}
\label{VZCtheta}
 V(Z(\mu))&\leq  \frac{V(T^n)}{\pi^{\frac{n+1}2}}
\int_C\theta(y)^{-1} e^{-\|\Theta(y)\|^2}\det\left( d\Theta(y)\right)\,dy\\
\label{VZ}
&\leq \frac{V(T^n)}{\pi^{\frac{n+1}2}}\int_{\R^{n+1}}
e^{-\|z\|^2}\,dz=\frac{n^{n/2}(n+1)^{(n+1)/2}}{n!}=V(T^n).
\end{align}

According to Lemma~\ref{ciuibig},
used for $v_i=\sqrt{\tilde{c}_i}\ \tilde{u}_i$, $i=1,\ldots,k$, we may assume that
\begin{equation}
\label{1ton}
\tilde{c}_1\cdots\tilde{c}_{n+1}
\det[\tilde{u}_1,\ldots,\tilde{u}_{n+1}]^2\geq {k \choose n+1}^{-1}.
\end{equation}
Then, in particular, the vectors $\tilde{u}_1,\ldots,\tilde{u}_{n+1}$ are linearly independent.
Since each factor on the left-hand side of (\ref{1ton}) is at most 1 (compare (\ref{cismall})),
the product of the remaining factors is at least ${k \choose n+1}^{-1}$. For Lemma~\ref{case12}, we define
\begin{equation}
\label{tauepsilon}
\varepsilon:=n^{60n}\tau^{1/4}<1 \qquad\mbox{and}\qquad
\omega:=\frac1{ 3^5 n^5 4^{n+1} n^{2n}}.
\end{equation}

In the following lemma, we adopt the assumptions and the notation from above.

\begin{lemma}
\label{case12}
Let the assumptions of Proposition \ref{volZstab} be satisfied. 
If  $i\in\{1,\ldots,k\}$, then  $\tilde{c}_i\leq
\omega^2 \varepsilon^2$ or  $\angle(\tilde{u}_i,\tilde{u}_j)\leq
\omega \varepsilon$ for some  $j\in \{1,\ldots, n+1\}$.
\end{lemma}
\proof If $i\in\{1,\ldots,n+1\}$, we can choose $j=i$ and then have $\angle(\tilde{u}_i,\tilde{u}_i)=0$. 
Thus it remains to consider the cases where $i\in\{n+2,\ldots,k\}$. For this, we proceed by contradiction and hence assume that
there is some  $i\in\{n+2,\ldots,k\}$ such that $\tilde{c}_i> \omega^2 \varepsilon^2$ and
$\angle(\tilde{u}_i,\tilde{u}_j)> \omega \varepsilon$ for all $j\in \{1,\ldots, n+1\}$.
Under this assumption, we will identify a subset $\Xi$ of $C$ with reasonably large volume such that
\begin{equation}
\label{thetayXi}
\theta(y)\geq 1+\gamma_0\,\varepsilon^4
\mbox{ \ for $y\in\Xi$},
\end{equation}
where  $\gamma_0:=n^{-18n-78}$ depends on $n$ (see (\ref{thetalow})). From this we will  then deduce a contradiction.

Since $\tilde{u}_1,\ldots,\tilde{u}_{n+1}$ are linearly independent,  there are uniquely determined
$\lambda_1,\ldots,\lambda_{n+1}\in\R$ such that
\begin{equation}\label{linind}
\tilde{u}_i=\lambda_1\tilde{u}_1+\ldots+\lambda_{n+1}\tilde{u}_{n+1}.
\end{equation}
We adjust the indices of
$\tilde{u}_1,\ldots,\tilde{u}_{n+1}$ so that
$$
\lambda_1\geq \ldots\geq \lambda_{n+1}.
$$
Since $\langle \tilde{u}_j,e_{n+1}\rangle=1/{\sqrt{n+1}}$ for
$j=1,\ldots,k$, we have $\lambda_1+\ldots+\lambda_{n+1}=1$, and thus we obtain 
 $\lambda_1\geq \frac1{n+1}$. Combining $\tilde c_1\le 1$, (\ref{1ton}), and \eqref{linind}, we thus conclude that
\begin{equation}
\label{2ni}
\tilde{c}_2\ldots\tilde{c}_{n+1}\tilde{c}_i
\det[\tilde{u}_2,\ldots,\tilde{u}_{n+1},\tilde{u}_i]^2\geq
\frac{\omega^2 \varepsilon^2}{(n+1)^2}{k \choose n+1}^{-1}
>\omega_0\varepsilon^2,
\end{equation}
where we define $\omega_0:=n^{-10n-30}$. The inequality on the right-hand side is 
confirmed by an elementary calculation, which is based on $k\le N+1$ and $n!\ge \sqrt{2\pi n}\ (n/e)^n$.

Next we construct the set $\Xi$ for which (\ref{thetayXi}) is satisfied.
The open convex cone
$$
C_0:=\left\{y\in\R^{n+1}:\,\langle y,e_{n+1}\rangle>\|y\|\ \frac{n}{\sqrt{n^2+1}}
\right\}
$$
satisfies $C_0\subset C$. In fact, if $y=x+re_{n+1}\in C_0$ with $x\in\R^n$ and $r>0$, then 
$$
r>\sqrt{\|x\|^2+r^2}\ \frac{n}{\sqrt{n^2+1}}.
$$
But this is equivalent to $\|x\|<r/n$, which in turn implies that $\langle x,u_i\rangle<r/\sqrt{n}$ for $i=1,\ldots,k$, hence $y\in C$. 

Writing $\alpha$ and $\beta$ to denote the acute angles with
 $\cos\alpha=\langle \tilde{u}_j,e_{n+1}\rangle=\frac1{\sqrt{n+1}}$,
$j=1,\ldots,k$, and $\cos\beta=\frac{n}{\sqrt{n^2+1}}$,
we have $\alpha-\beta<\angle(y,\tilde{u}_j)<\alpha+\beta$
for $y\in C_0$ and $j=1,\ldots,k$. For $y\in C_0$ and $j=1,\ldots,k$,
we deduce that
\begin{align}
\langle y,\tilde{u}_j\rangle&<\|y\|\ \frac{n+\sqrt{n}}{\sqrt{(n^2+1)(n+1)}}
<\|y\|\ \frac2{\sqrt{n}},\label{ineqa}\\
\langle y,\tilde{u}_j\rangle&>\|y\|\ \frac{n-\sqrt{n}}{\sqrt{(n^2+1)(n+1)}}
>\|y\|\  \frac1{5\sqrt{n}}.\label{ineqb}
\end{align}
To verify the left inequality in \eqref{ineqa}, we consider $y=x+re_{n+1}\in C_0$ with $\|y\|=1$. 
Then $\|x\|^2+r^2=1$ and $r>n/\sqrt{n^2+1}$. Hence
$$
\langle y,\tilde u_j\rangle =-\sqrt{\frac{n}{n+1}}\langle x,u_j\rangle +\frac{r}{\sqrt{n+1}}
\le \sqrt{\frac{n}{n+1}}\sqrt{1-r^2} +\frac{r}{\sqrt{n+1}}=:f(r).
$$
Since $f$ is decreasing for $r\ge n/\sqrt{n^2+1}$, the assertion follows. Similarly, 
$$
\langle y,\tilde u_j\rangle \ge -\sqrt{\frac{n}{n+1}}\sqrt{1-r^2} +\frac{r}{\sqrt{n+1}}=:g(r)
$$
and $g$ is increasing for $r\ge n/\sqrt{n^2+1}$, which yields the first inequality in \eqref{ineqb}.

We also observe that the section
$\{y\in C_0:\,\langle y,e_{n+1}\rangle=t\}$ is an $(n-1)$-ball of radius $t/n$
for $t>0$. Now we are ready to define
$$
\Xi:=\left\{y\in C_0:\,20\sqrt{n}<\langle y,e_{n+1}\rangle<40\sqrt{n}
\mbox{ and }
\langle y,\tilde{u}_i-\tilde{u}_1\rangle>\frac{\omega\varepsilon}{\sqrt{n}}\right\}.
$$
Since by assumption $\|\tilde{u}_i-\tilde{u}_1\|>\omega\varepsilon/2$,
$\Xi$ contains a right cylinder of height $20\sqrt{n}$ whose base is 
an $(n-1)$-dimensional regular simplex $S_*$ of circumradius $1/\sqrt{n}$.
Let $S_0$ be an $n$-dimensional regular simplex  whose facet is $S_*$. Since
the height of $S_0$ is less than $2/\sqrt{n}$, we have
\begin{equation}
\label{Xivol}
V(\Xi)>\frac{n\ 20\sqrt{n}}{2/\sqrt{n}}\,V(S_0)= 
\frac{10\ n^2}{n^{3n/2}}\,V(T^n).
\end{equation}
Using \eqref{ineqa} and \eqref{ineqb}, we also get
\begin{equation}
\label{Xiy}
4<\langle y,\tilde{u}_j\rangle<120
\mbox{ \ for $y\in \Xi$ and $j=1,\ldots,k$.}
\end{equation}

For $y\in \Xi$, we estimate $\theta(y)$ from below using
 the $n$-tuples $(1,\ldots,n+1)$ and $(2,\ldots,n+1,i)$
of indices in (\ref{thetay}) (note that in addition to \eqref{2ni} we also have $\binom{k}{n+1}^{-1}\ge \omega_0$). We deduce by first
applying (\ref{1ton}), (\ref{2ni}) and  Lemma~\ref{xab}, secondly
$\varphi'(\langle y,\tilde{u}_j\rangle)<1$ for $j=1,\ldots,k$
(see Lemma~\ref{phi}), and thirdly
by $\langle y,\tilde{u}_i-\tilde{u}_1\rangle>\frac{\omega\varepsilon}{\sqrt{n}}$
and $\varphi''(t)<-12^{-4}$ for $4<t<120$ (see Lemma~\ref{phi}) that
\begin{align}
\nonumber
\theta(y)&\geq 1+\frac12\ 
\frac{(\varphi'(\langle y,\tilde{u}_1\rangle)-\varphi'(\langle y,\tilde{u}_i\rangle))^2}
{2(\varphi'(\langle y,\tilde{u}_1\rangle)+\varphi'(\langle y,\tilde{u}_i\rangle))^2}
\ \omega_0\varepsilon^2\\
\nonumber
&> 1+\frac{(\varphi'(\langle y,\tilde{u}_1\rangle)-\varphi'(\langle y,\tilde{u}_i\rangle))^2}{16}
\  \omega_0\varepsilon^2\\
\label{thetalow}
&> 1+\frac{\omega^2\omega_0}{16\ n\ 12^8}\ \varepsilon^4
>1+n^{-18n-78}\ \varepsilon^4.
\end{align}

According to (\ref{Xiy}) and  Lemma~\ref{phi}, if $y\in \Xi$ and $j=1,\ldots,k$, then
$\varphi(\langle y,\tilde{u}_j\rangle)^2<120$ and
$\varphi'(\langle y,\tilde{u}_j\rangle)>\frac{1}{33}$.
It follows from (\ref{firstterm}) and (\ref{secondterm}), taking into account (\ref{John3n+1}), that
\begin{align}
\nonumber
e^{-\|\Theta(y)\|^2}\det \left(d\Theta(y)\right)&\geq
\exp\left(-\sum_{j=1}^k \tilde{c}_j\varphi(\langle y,\tilde{u}_j\rangle)^2  \right)
\prod_{j=1}^k\varphi'(\langle y,\tilde{u}_j\rangle)^{\tilde{c}_j}\\
\label{Tlow}
&\ge  e^{-120(n+1)}\, {33^{-(n+1)}}\ge e^{-124(n+1)}\geq e^{-186n}.
\end{align}
Recall that $\gamma_0=n^{-18n-78}$ and observe that (\ref{thetalow}) implies that
\begin{equation}\label{eqc}
1-\theta(y)^{-1}\ge \frac{\gamma_0\varepsilon^4}{1+\gamma_0\varepsilon^4}\ge \frac{1}{2}\gamma_0\varepsilon^4.
\end{equation}
Now we use  (\ref{Xivol}), (\ref{Tlow}) and \eqref{eqc}, and argue as for 
 (\ref{VZCtheta}) and (\ref{VZ}), to obtain
\begin{align*}
 V(Z(\mu))
&\le \frac{V(T^n)}{\pi^{\frac{n+1}2}}
\int_C  e^{-\|\Theta(y)\|^2}\det\left( d\Theta(y)\right)\,dy\\
&\qquad  -\frac{V(T^n)}{\pi^{\frac{n+1}2}}
\int_C\left(1-\theta(y)^{-1}\right) e^{-\|\Theta(y)\|^2}\det\left( d\Theta(y)\right)\,dy\\
&\le V(T^n)-\frac{V(T^n)}{\pi^{\frac{n+1}2}}
\int_\Xi\left(1-\theta(y)^{-1}\right) e^{-\|\Theta(y)\|^2}\det\left( d\Theta(y)\right)\,dy\allowdisplaybreaks\\
&\le V(T^n)\left[1-\frac{1}{\pi^{\frac{n+1}2}}
\int_\Xi \tfrac{1}{2}\gamma_0 \varepsilon^4 e^{-186 n}\,dy\right]\allowdisplaybreaks\\
&\le V(T^n)\left[1-\frac{V(\Xi)}{2\pi^{\frac{n+1}2}}\gamma_0 \varepsilon^4 e^{-186 n}\right]\allowdisplaybreaks\\
&\le V(T^n)\left[1-\frac{5n^2V(T^n)}{n^{\frac{3n}{2}}\pi^{\frac{n+1}2}}\gamma_0 \varepsilon^4 e^{-186 n}\right]\allowdisplaybreaks\\
&\le \left(1-n^{-240n}\varepsilon^4\right) V(T^n)=(1-\tau)V(T^n),
\end{align*}
where we used  (\ref{tauepsilon}) in the last step.  
This contradicts the assumptions of Proposition \ref{volZstab}, and hence
 proves Lemma~\ref{case12}.
\proofbox

\noindent{\bf Proof of Proposition~\ref{volZstab}: } 
For $i=1,\ldots,k$, we define $\tilde v_i:=\sqrt{\tilde c_i}\tilde u_i\in\R^{n+1}$, hence $\|\tilde v_i\|=\sqrt{\tilde c_i}$. Lemma \ref{case12} 
ensures that the assumptions for the application of Lemma \ref{almostort} are satisfied for $\tilde v_1,\ldots,\tilde v_k$ in $\R^{n+1}$ with 
$\eta=\omega\varepsilon <1/(3\sqrt{k})$. Hence, by Lemma \ref{almostort} there is an orthonormal basis $\bar w_1,\ldots,\bar w_{n+1}$ of $\R^{n+1}$ 
such that $\angle(\tilde v_i,\bar w_i)<3\sqrt{k}\omega\varepsilon$ for $i=1,\ldots,n+1$. Writing $\alpha_i=\angle (e_{n+1},\bar w_i)$ 
and  $\beta_i=\angle (e_{n+1},\tilde v_i)=\angle (e_{n+1},\tilde u_i)$, we get 
$$
\left|\langle e_{n+1},\bar w_i\rangle - \frac{1}{\sqrt{n+1}}\right|=|\cos\alpha_i-\cos\beta_i|\le |\alpha_i-\beta_i|\le 
\angle (\bar w_i,\tilde v_i)<3\sqrt{k}\omega\varepsilon. 
$$
Since $3\sqrt{k}\omega\varepsilon<1/(2(n+1))$, we can apply Lemma \ref{basis}, which yields the existence of an orthonormal basis 
$\tilde w_1,\ldots,\tilde w_{n+1}$ in $\R^{n+1}$ such that $\langle e_{n+1},\tilde w_i\rangle=1/\sqrt{n+1}$ 
and $\angle(\tilde w_i,\bar w_i)\le (n+1)3\sqrt{k}\omega\varepsilon$. But then 
$$
\angle(\tilde w_i,\tilde u_i)\le \angle(\tilde w_i,\bar w_i)+\angle(\bar w_i,\tilde u_i,)
\le 3(n+1)\sqrt{k}\omega\varepsilon+3\sqrt{k}\omega\varepsilon\le 8n^2\omega\varepsilon.
$$ 

For $i=1,\ldots,n+1$, we define 
$$
w_i=\sqrt{\frac{n+1}n}\,\left(-\tilde{w}_i+\sqrt{\frac1{n+1}}\,e_{n+1}  \right)\in \R^n,
$$
and hence there exists a regular simplex $S$ whose facets touch $B^n$
at $w_1,\ldots,w_{n+1}$. Subsequently, we use that
$$
1-\frac12\,t^2<\cos t<1-\frac{3}{8}\,t^2
\mbox{ \ for $t\in(0,1)$}.
$$
Since
$$
1-\langle w_i,u_i\rangle=\frac{n+1}n\,(1-\langle \tilde{w}_i, \tilde{u}_i\rangle)
\leq \frac{n+1}n\ \frac12(8n^2\,\omega \varepsilon)^2\leq
48 n^4\omega^2\varepsilon^2,
$$
we deduce that $\angle( w_i,u_i)<12 n^2\omega\varepsilon$ for $i=1,\ldots,n+1$.

We observe that  $\gamma=9\cdot 2^{n+2}n^{2n+2}$ from  Lemma~\ref{closefarsimplex}
and $\omega=(3^5 n^5 4^{n+1}n^{2n})^{-1}$ satisfy
\begin{equation}
\label{finallyclose0}
\frac{1}{9 n2^{n-2}}\le 12 \gamma n^2\omega\le \frac{1}{9n},
\end{equation}
and claim that
\begin{equation}
\label{finallyclose}
\delta_H({\rm supp}\,\mu,{\rm supp}\,\mu_S)< 12 \gamma n^2\omega\varepsilon\le \frac{1}{9n}\varepsilon = \frac{1}{9n} n^{60n}\tau^{1/4}.
\end{equation}
Let us suppose that contrary to (\ref{finallyclose}), there exists some $i\in\{n+2,\ldots,k\}$ such that
$\angle (u_i,w_j)\ge  12 \gamma n^2\omega\varepsilon$ for  $j=1,\ldots,n+1$.
To apply Lemma~\ref{closefarsimplex}, we note that $\varepsilon<1$ and (\ref{finallyclose0})
 yield that
$12 n^2\omega\varepsilon<\gamma^{-1}$. Since
$\varepsilon=n^{60n}\tau^{1/4}>n^{240n}\tau$, we conclude from (\ref{finallyclose0}) that
$$
V(Z(\mu))\leq \left(1-\frac{12 \gamma n^2\omega\varepsilon }{2^{n+2}n^{2n}}\right)V(T^n)<
(1-\tau)V(T^n).
$$
This contradicts the condition on $\mu$, and hence implies (\ref{finallyclose}).
Finally, combining (\ref{finallyclose}) and Lemma~\ref{closesimplex} yields
$d(Z(\mu))<n^{60n}\tau^{1/4}$.\proofbox

\section{Proofs of Theorems~\ref{vol} and \ref{BM}}
\label{secisoperimetric}

We assume that $B^n$ is the ellipsoid of maximal volume inside the convex body $K$ in $\R^n$,
and hence there exist
$u_1,\ldots,u_k\in S^{n-1}\cap \partial K$ and $c_1,\ldots,c_k>0$ such that
$\sum_{i=1}^kc_i u_i\otimes u_i={\rm Id}_n$ and
$\sum_{i=1}^kc_i u_i=o$, where
\begin{equation}
\label{Johnbound0}
n+1\leq k\leq n(n+3)/2.
\end{equation}
We write $Z$ to denote the circumscribed polytope whose faces touch $B^n$ at
$u_1,\ldots,u_k$; namely,
$$
Z=\{x\in\R^n:\,\langle x,u_i\rangle\leq 1,\;i=1,\ldots,k\}.
$$

For any $x\in\partial K$, let $u_x$  denote an exterior
unit normal at $x$, which is unique (almost everywhere) and measurable with respect
to the $(n-1)$-dimensional Hausdorff-measure on $\partial K$. We note that
\begin{equation}
\label{SVint}
V(K)=\int_{\partial K}\frac{\langle x,u_x\rangle}n\,dx\geq \frac{S(K)}n.
\end{equation}
It follows from (\ref{SVint}) that
\begin{equation}
\label{SVbasic}
\frac{S(K)^n}{V(K)^{n-1}}\leq n^nV(K)\leq n^n V(Z)
\leq n^nV(T^n)=\frac{S(T^n)^n}{V(T^n)^{n-1}}.
\end{equation}

\begin{lemma}
\label{ZK}
Let  $\varepsilon\in(0,1)$. 
\begin{description}
\item{{\rm (i)}} If $d(Z)\leq \varepsilon/(4n^2)$ and
 $\delta_{\rm BM}(K,T^n)\geq \varepsilon$, then
$$
\frac{S(K)^n}{V(K)^{n-1}}\le \left(1-\frac{1}{e^2}\left(\frac{\varepsilon}{e}\right)^n\right)
\frac{S(T^n)^n}{V(T^n)^{n-1}}.
$$
\item{{\rm (ii)}} If $d(Z)\leq \varepsilon/(4n^2)$ and
$\delta_{\rm vol}(K,T^n)\geq \varepsilon$, then
$$
\frac{S(K)^n}{V(K)^{n-1}}\le \left(1-\frac{\varepsilon}{8}\right)
\frac{S(T^n)^n}{V(T^n)^{n-1}}.
$$
\end{description}
\end{lemma}
\proof Let $\gamma:=1/(4n^2)$. Then we  may assume that
\begin{equation}
\label{ZTncond}
e^{-\gamma\varepsilon}\,T^n\subset Z\subset e^{\gamma\varepsilon}\,T^n.
\end{equation}
Hence, we have $4n\gamma\varepsilon\le\frac1{n}$.

For the proof of (ii), we first choose $\lambda>0$ such that $V(T^n)=V(\lambda Z)$. 
Then (\ref{ZTncond}) yields that 
$e^{-\gamma\varepsilon}\leq \lambda\leq e^{\gamma\varepsilon}$. Therefore, again by (\ref{ZTncond}) we obtain
\begin{align*}
\delta_{\rm vol}(Z,T^n)&\leq \frac{V((\lambda Z)\Delta T^n)}{V(T^n)}
\leq
\lambda^n e^{n\gamma\varepsilon}-
\lambda^n e^{-n\gamma\varepsilon}\\
&\le 2n\gamma\varepsilon \lambda^n e^{n\gamma\varepsilon} \le 2n\gamma\varepsilon   e^{2n\gamma\varepsilon}\le 2n\gamma\varepsilon (1+4n\gamma\varepsilon)\\
&\leq 4n\gamma\varepsilon\le\varepsilon/2,
\end{align*}
where we used that $e^t\le 1+2t$ for $0\le t\le 1/2$. 

Let $\eta,\nu\geq 0$ satisfy $V(K)=V(\eta Z)$ and $V(Z)=(1+\nu)V(K)$,
and hence $\eta=(1+\nu)^{-1/n}$.
 It follows from $\delta_{\rm vol}(K,T^n)\geq \varepsilon$ that
$$
\varepsilon/2\leq
\delta_{\rm vol}(Z,K)\leq \frac{V((\eta Z)\Delta K)}{V(K)}\leq
\frac{2V(Z\backslash K)}{V(K)}
\leq  2\nu,
$$
and hence (\ref{SVbasic}) yields that
\begin{align*}
\frac{S(K)^n}{V(K)^{n-1}}&\leq n^n V(K)=n^n (1+\nu)^{-1}V(Z)\le (1+\nu)^{-1}\frac{S(T^n)^n}{V(T^n)^{n-1}} \\
&\le \left(1+\frac{\varepsilon}{4}\right)^{-1}  
\frac{S(T^n)^n}{V(T^n)^{n-1}}\leq \left(1-\frac{\varepsilon}{8}\right) \frac{S(T^n)^n}{V(T^n)^{n-1}}.
\end{align*}

We turn to (i).
It follows from $\delta_{\rm BM}(K,T^n)\ge \varepsilon$ and (\ref{ZTncond})
that there is a vertex $v$ of $T^n$ such that
$$
e^{\gamma\varepsilon-\varepsilon}v\not\in {\rm int}\,K.
$$
In particular, there exists a half-space $H^+$ containing $e^{\gamma\varepsilon-\varepsilon}v$, and disjoint from ${\rm int}\,K$.
Since $p=e^{\gamma\varepsilon-\varepsilon}v$ is the centroid of the simplex
$p+\lambda T^n\subset e^{-\gamma\varepsilon}T^n$ for
$\lambda:=e^{-\gamma\varepsilon}-e^{\gamma\varepsilon-\varepsilon}$, a result by B.~Gr\"unbaum \cite[p.~1260, (iii)]{Gru60} yields that
$$
V(H^+\cap(p+\lambda T^n))>\frac{\lambda^n}e\, V(T^n).
$$
Therefore, using  (\ref{ZTncond})  we deduce  that
\begin{align*}
V(Z\setminus K)&\geq V(H^+\cap(e^{-\gamma\varepsilon}T^n))>\frac{\lambda^n}e\,V(T^n)
=\frac{e^{-n\gamma \varepsilon}(1-e^{2\gamma\varepsilon-\varepsilon})^n}e V(T^n)\\
&\ge \frac{1}{e^2}\left(\frac{\varepsilon}{e}\right)^n V(T^n).
\end{align*}
Hence, by (\ref{SVbasic}) we get
$$
V(K)+\frac{1}{e^2}\left(\frac{\varepsilon}{e}\right)^n V(T^n)\le V(Z)\le V(T^n),
$$
and therefore
$$
V(K)\le \left(1-\frac{1}{e^2}\left(\frac{\varepsilon}{e}\right)^n \right)V(T^n).
$$
Now the proof can be completed as in the previous case by using once again (\ref{SVbasic}). 
\mbox{ \ }\proofbox

\noindent{\em Proofs of Theorems~\ref{vol} and \ref{BM}: }
If $d(Z)>\varepsilon/(4n^2)$, then Proposition~\ref{volZstab} can be applied by (\ref{Johnbound0}),
and implies that
$$
V(Z)\leq(1-4^{-4}n^{-248n}\varepsilon^4)V(T^n)\leq (1-n^{-250n}\varepsilon^4)V(T^n).
$$
 In turn, we conclude
Theorem~\ref{BM} and  Theorem~\ref{vol} by (\ref{SVbasic}).

If  $d(Z)\leq \varepsilon/(4n^2)$, then Lemma~\ref{ZK} (i) yields Theorem~\ref{BM},
and  Lemma~\ref{ZK} (ii) implies Theorem~\ref{vol}.
 \proofbox
 
For the sake of completeness we provide the following fact, which is mentioned in the introduction.
 
\begin{lemma}
Let $K,M$ be convex bodies in $\R^n$. Then $\delta_{\rm vol}(K,M)\le 2 e^{n^2}\delta_{\rm BM}(K,M)$ and $\delta_{\rm BM}(K,M)\le \gamma\ \delta_{\rm vol}(K,M)^{\frac{1}{n}}$, 
where $\gamma$ is a constant which depends on $n$. 
\end{lemma}

\proof The assertions follow from  \cite[Section 5]{BHenk}. Since the first assertion is used explicitly (in the introduction) and  the definitions of the distances used here differ from those given in \cite{BHenk}, we outline the short argument for the first inequality.

Since $\delta_{\rm vol}$ and $\delta_{\rm BM}$ are translation invariant in both arguments, we can assume that $0\in K,M$ and $K\subset M\subset e^\delta K$, where  $\delta:=\delta_{\rm BM}$, and therefore $V(K)\le V(M)\le e^{n\delta} V(K)$ or 
$$
1\le \left(\frac{V(M)}{V(K)}\right)^{\frac{1}{n}}\le e^\delta.
$$
Thus we conclude that 
$$
e^{-\delta} K_0\subset M_0\subset e^\delta K_0,
$$
where $K_0:=V(K)^{-\frac{1}{n}} K$ and $M_0:=V(M)^{-\frac{1}{n}} M$. But then
$$
V(K_0\Delta M_0)\le V((e^\delta K_0)\setminus K_0)+V((e^\delta M_0)\setminus M_0)\le 2\left(e^\delta-1\right)\le 2\delta e^\delta.
$$
Now the assertion follows since $\delta_{\rm BM}(K,M)\le n^2$. 
\proofbox

\section{Proof of Theorem \ref{2D}}
\label{sec2DProof}

Throughout the proof, we have $n=2$. The argument is based on \cite{Gustin1953}, which we briefly recall. For a convex body $K$ in 
$\R^n$ and $u\in S^{n-1}$, we write $H^-(K,u)$ for the supporting half-space of $K$  which contains $K$ and has exterior unit normal $u$, and $H(K,u)$ for 
its bounding hyperplane.

For the proof, we assume that 
\begin{equation}\label{reverse2d}
\text{ir}(K)\ge (1-\varepsilon)\text{ir}(T^2).
\end{equation} 
Let $\text{IR}(K):=S(K)^2/V(K)$ for a convex body $K$ in $\R^2$.  Then $\text{ir}(T^2)=\text{IR}(T^2)$. 
 Let $T_1$ be a triangle of maximal area contained in $K$. We can assume that $T_1$ is a regular triangle centred at $0$ with height $1$, whose vertices are denoted by 
 $p_1,p_2,p_3$. Let $u_1,u_2,u_3\in S^1$ 
denote the exterior normal vectors of the edges of $T_1$. Then the lines $H(T_1,-u_i)$, $i=1,2,3$, pass through the vertices of $T_1$ and bound a regular triangle $T_2$ of height $2$ which contains $K$. 
Choose 
$q_i\in K\cap H(K,u_i)$ and let $x_i\in [0,1]$ be the distance of $q_i$ from $H(T_1,u_i)$ for $i=1,2,3$. Then   
$$
T_1\subset P_1:=\text{conv}\{p_1,p_2,p_3,q_1,q_2,q_3\}\subset K\subset \bigcap_{i=1}^3H^-(K,u_i)\cap \bigcap_{i=1}^3H^-(K,-u_i)=:P_2\subset T_2.
$$
Let $x:=(x_1+x_2+x_3)/3\in [0,1]$. Elementary geometric arguments show (see \cite{Gustin1953}) that 
$$
S(P_2)=(1+x)S(T_2)\qquad\text{and}\qquad V(P_1)=(1+3x)V(T_2),
$$
and therefore
$$
\text{ir}(K)\le \text{IR}(K)\le \frac{S(P_2)^2}{V(P_1)}\le \left(1-\frac{x(1-x)}{1+3x}\right)\text{ir}(T^2).
$$
From \eqref{reverse2d} we conclude that $(1+3x)^{-1}x(1-x)\le \varepsilon$, and thus $x(1-x)\le 4\varepsilon$. 

If $x\le 1/2$, then $x\le 8\varepsilon$ and thus $x_i\le 24\varepsilon$ for $i=1,2,3$. If $x\ge 1/2$, then in fact $x\ge 1-8\varepsilon$ 
and hence $x_i\ge 1-24\varepsilon$ for $i=1,2,3$. In the first case, we conclude that 
$$
T_1\subset K\subset P_2\subset(1+72\varepsilon)T_1,
$$
which implies
$$
\delta_{\text{BM}}(K,T^2)\le \ln(1+72\varepsilon)\le 72\varepsilon. 
$$
In the second case, we find a regular triangle $T$ centred at $0$ and homothetic to $T_2$ such that $T\subset K\subset T_2$ whose edges have distance at 
least $(2/3)-24(2/3)\sqrt{3}\ \varepsilon$ from $0$. This shows that 
$$
\delta_{\text{BM}}(K,T^2)\le \ln\left(\frac{1}{1-24\sqrt{3}\varepsilon}\right)\le 72\varepsilon 
$$
for $\varepsilon\le 1/72$. This completes the proof in both cases.

\section{Isotropic measures: proof of Theorem \ref{Zmustab}}
\label{secisotropic}

Our proof of  Theorem~\ref{Zmustab} will be based on Proposition \ref{volZstab}. For this reason we have to ensure that we can switch from a centred, isotropic measure 
$\mu$ on $S^{n-1}$ to a discrete, centred, isotropic measure on $S^{n-1}$ with support contained in the support of $\mu$ and whose support has bounded cardinality. That this 
can indeed be achieved is shown by the following lemma. 

Recall that $N=n(n+3)/2$.

\begin{lemma}\label{cara}
Let $\mu$ be a centred, isotropic measure on $S^{n-1}$. Then there exists a discrete, centred,  
isotropic measure $\mu_0$ on $S^{n-1}$ such that  $\text{\rm supp}\ \mu_0\subset \text{\rm supp}\ 
\mu$ and the cardinality of $\text{\rm supp}\ \mu_0$ is at most $N+1$.
\end{lemma}

\proof We consider the map $F:\text{\rm supp}\ \mu\to \R^N$ given by $F(u):=(u\otimes u,u)$. 
Here we interpret $u\otimes u$ as the upper triangular part (including the main diagonal) of the 
symmetric matrix $u\otimes u$, and thus we identify the vectors $(u\otimes u,u)$ with vectors in $\R^N$. 
Since $\text{\rm supp} \ \mu\subset S^{n-1}$ is compact and $F$ is continuous, 
the image set $F(\text{\rm supp}\ \mu)\subset\R^N$ is compact as well. Then also the convex hull of this image set, 
$\text{conv}(F(\text{\rm supp}\ \mu))\subset\R^N$ is compact. The probability measure $\bar\mu:=\mu/n$ has the same support 
as $\mu$ and satisfies
$$
\left(\int_{S^{n-1}}u\otimes u\,d\bar\mu(u),\int_{S^{n-1}}u\, d\bar\mu(u)\right)=\left(\frac{1}{n}\text{Id}_n,0\right)\in\R^N.
$$
Let $\mathcal{D}_l$ be a decomposition of $S^{n-1}$ into finitely many disjoint Borel sets of diameter at most $1/l$, $l\in\N$. 
We put $\mathcal{D}_l^*:=\{\Delta\in \mathcal{D}_l: \Delta\cap\text{supp}\ \bar\mu\neq\emptyset\}$. For $\Delta\in \mathcal{D}_l^*$, 
we fix some $v_\Delta\in \Delta\cap\text{supp}\ \bar\mu$. Then 
$$
\bar\mu_l:=\sum_{\Delta\in\mathcal{D}_l^*}\bar\mu(\Delta)\delta[v_\Delta]
$$
is a discrete probability measure on $S^{n-1}$ and $\text{supp}\ \bar\mu_l\subset \text{supp}\ \bar\mu$. Moreover, $\bar\mu_l\to \bar\mu$ 
in the weak topology as $l\to\infty$. Therefore, we conclude that
$$
\sum_{\Delta\in\mathcal{D}_l^*}\bar\mu(\Delta)\left(v_\Delta\otimes v_\Delta,v_\Delta\right)=
\left(\int_{S^{n-1}}v\otimes v\,d\bar\mu_l(v),\int_{S^{n-1}}v\, d\bar\mu_l(v)\right)
\to\left(\frac{1}{n}\text{Id}_n,0\right) 
$$
in $\R^N$ as $l\to\infty$. This shows that 
$$
\left(\frac{1}{n}\text{Id}_n,0\right)\in \text{cl}\ \text{conv}(F(\text{supp}\ \bar\mu))=\text{conv}(F(\text{supp}\ \bar\mu)).
$$ 
By Carath\'eodory's theorem (see, e.g., \cite[Theorem 1.1.4]{Sch14}) there exist $k\le N+1$ vectors $u_1,\ldots,u_k\in\text{supp}\ \bar\mu\subset S^{n-1}$ such that 
$$
\left(\frac{1}{n}\text{Id}_n,0\right)\in \text{conv}(F(\{u_1,\ldots,u_k\})),
$$
that is, there exist $\alpha_1,\ldots,\alpha_k\ge 0$ with $\alpha_1+\ldots+\alpha_k=1$ such that 
$$
\left(\frac{1}{n}\text{Id}_n,0\right)=\sum_{i=1}^k\alpha_iF(u_i)=\sum_{i=1}^k \alpha_i (u_i\otimes u_i,u_i).
$$
This shows that with $c_i:=n\alpha_i$ for $i=1,\ldots,k$ the measure
$$
\mu_0:=\sum_{i=1}^kc_i\delta[u_i]
$$
satisfies all requirements.
\proofbox

\bigskip

For the proof of Theorem \ref{Zmustab} we can assume that $\varepsilon\in (0,n^{-268n})$, since otherwise $n^{70n}\varepsilon^{\frac{1}{4}}\ge n^{3n}$ and the assertion is trivial. For the given measure $\mu$ there is a measure $\mu_0$ as described in Lemma  \ref{cara}. Combined with the assumption of Theorem \ref{Zmustab} this yields that 
$$
(1-\varepsilon)V(T^n)\le V(Z(\mu))\le V(Z(\mu_0)).
$$
Hence we can apply Proposition \ref{volZstab} and obtain a regular simplex $S$ circumscribed about $B^n$ with contact points $w_1,\ldots,w_{n+1}$ and such that
\begin{equation}\label{closediscrete}
\delta_H(\text{supp}\ \mu_0,\text{supp}\ \mu_S)\le n^{60n}\varepsilon^{\frac{1}{4}}.
\end{equation}
If $\text{supp}\ \mu_0=\text{supp}\ \bar\mu$, the proof is finished. Hence, let $u^*\in \text{supp}(\bar\mu)\setminus \text{supp}(\mu_0)$ and let $Z^*$ be the polytope circumscribed to $B^n$ with contact points $\text{supp}(\mu_0)\cup\{u^*\}$. Then we have
$$
(1-\varepsilon)V(T^n)\le V(Z(\mu))\le V(Z^*).
$$
Let $\eta:=n^{60n}\varepsilon^{\frac{1}{4}}<\gamma^{-1}=(9\cdot 2^{n+2}n^{2n+2})^{-1}$. From \eqref{closediscrete} we conclude that we can assume that 
$\text{supp}\ \mu_0=\{u_1,\ldots,u_{k}\}$, $k\ge n+1$, with $\angle(u_i,w_i)\le \eta$ for $i=1,\ldots,n+1$. Assume that 
$\angle (u^*,w_i)\ge \gamma\eta$ for $i=1,\ldots,n+1$. Then Lemma \ref{closefarsimplex} implies that 
$$
(1-\varepsilon)V(T^n)\le V(Z^*)\le \left(1-\frac{\gamma\eta}{2^{n+2}n^{2n}}\right)V(T^n),
$$
and therefore $\gamma\eta\le 2^{n+2}n^{2n}\varepsilon$, which contradicts $\varepsilon\le 1$. This 
shows that $\angle (u^*,w_i)<\gamma\eta$ for some $i\in\{1,\ldots,n+1\}$. Since $\gamma\eta\le n^{67n}\varepsilon^{\frac{1}{4}}$, it finally follows that  
$\delta_H(\text{supp}\ \bar\mu,\text{supp}\ \mu_S)\le n^{67n}\varepsilon^{\frac{1}{4}}$, which proves the theorem.
\proofbox

Finally, we justify the remark following Theorem \ref{Zmustab} by establishing the next lemma. For $w\in S^{n-1}$ and $\varepsilon\ge 0$, we consider 
 $U(w,\varepsilon):=\{u\in S^{n-1}:\angle (u,w)\le \varepsilon\}$, that is, the closed spherical (geodesic) ball with centre $w$ and radius $\varepsilon$.

\begin{lemma}\label{massbound}
Let $S$ be a regular simplex circumscribed about $B^n$ with contact points $w_1,\ldots,w_{n+1}\in S^{n-1}$, let $\mu$ be a centred, isotropic Borel measure 
on $S^{n-1}$, and let $\varepsilon\in(0,1/2)$. If $\delta_H(\text{\rm supp}\ \mu,\text{\rm supp}\ \mu_S)\le \varepsilon$, then
$$
\left|\mu(U(w_i,\varepsilon))-\frac{n}{n+1}\right|\le 2n\varepsilon,\qquad i=1,\ldots,n+1.
$$ 
\end{lemma}

\proof Let the map $G:S^{n-1}\to S^n$ be defined by
$$
G(u):=-\sqrt{\frac{n}{n+1}} \ u+\sqrt{\frac{1}{n+1}}\ e_{n+1}.
$$
Since $\mu$ is centred and isotropic, we obtain 
$$
{\rm Id}_{n+1}=\frac{n+1}n\int_{S^{n-1}}G(u)\otimes G(u)\, d\mu(u). 
$$
By assumption,  $\text{supp}\ \mu \subset \bigcup_{i=1}^{n+1} U(w_i,\varepsilon)$ and the union is disjoint. For $u\in U(w_i,\varepsilon)$ and $x\in S^{n}$, 
using the triangle and the Cauchy-Schwarz inequality as well as the fact that $G(u), G(w_i)$ and $x$ are unit vectors, we get
$$
\left\| \langle G(u),x\rangle G(u)-\langle G(w_i),x\rangle G(w_i)\right\|\le 2\|G(u)-G(w_i)\|\le 2\|u-w_i\|\le 2\varepsilon.
$$
Hence, for any $x\in S^n$,
\begin{align*}
&\left\|x-\frac{n+1}{n}\sum_{i=1}^{n+1}\mu(U(w_i,\varepsilon))\langle G(w_i),x\rangle G(w_i)\right\|\\
&\qquad = \frac{n+1}{n}\left\|\int_{S^{n-1}}\langle G(u),x\rangle G(u)\, d\mu(u)-\sum_{i=1}^{n+1}\mu(U(w_i,\varepsilon))\langle G(w_i),x\rangle G(w_i)
\right\|\allowdisplaybreaks\\
&\qquad\le  \frac{n+1}{n}\sum_{i=1}^{n+1}\int_{U(w_i,\varepsilon)}\left\| \langle G(u),x\rangle G(u)-\langle G(w_i),x\rangle G(w_i)
\right\|\, d\mu(u)\allowdisplaybreaks\\
&\qquad\le \frac{n+1}{n} 2\varepsilon \sum_{i=1}^{n+1}\mu(U(w_i,\varepsilon))=2(n+1)\varepsilon.
\end{align*}
The special choice $x=G(w_i)$, for some $i\in \{1,\ldots,n+1\}$, together with the fact that $G(w_1),\ldots,G(w_{n+1})$ is an orthonormal basis of $\R^{n+1}$ then yields 
$$
|1-((n+1)/n)\mu(U(w_i,
\varepsilon)|\le 2(n+1)\varepsilon,
$$ 
from which the assertion follows.
\proofbox

Let the assumptions of Lemma \ref{massbound} be satisfied. 
Furthermore, let $f:S^{n-1}\to\R$ be lipschitz with lipschitz constant $\|f\|_L$. Here the definition of the lipschitz constant 
is based on the geodesic distance on $S^{n-1}$.  Since $\mu$ and $\mu_S$ have the same total measure $n$, we 
can replace $f$ by $f-f(e_1)$ in the following estimation, and therefore we can assume that the sup norm $\|f\|_\infty$ of $f$ 
satisfies $\|f\|_\infty\le 4\|f\|_L$. Thus, we get
\begin{align*}
\left|\int_{S^{n-1}}f\, d\mu -\int_{S^{n-1}}f\, d\mu_S\right|& \le \sum_{i=1}^{n+1}\int_{U(w_i,\varepsilon)}|f-f(w_i)|\, d\mu +\sum_{i=1}^{n+1}|f(w_i)|2n\varepsilon\\
&\le \|f\|_L\varepsilon n+\|f\|_\infty 2n(n+1)\varepsilon\\
&\le 13n^2\varepsilon \|f\|_L,
\end{align*}
which yields the asserted bound for the Wasserstein distance $d_W(\mu,\mu_S)$.

\bigskip

\section{Proof of Theorem \ref{planariso}}
\label{2dcase}

We state the next lemma in general dimensions although we will need it only in the plane.

\begin{lemma}\label{lowerbound}
Let $\mu$ be a centred and isotropic Borel measure on $S^{n-1}$. Let $v\in S^{n-1}$ be given. Then there is some $u^*\in \text{\rm supp}\ \mu$ such that $
\langle u^*,v\rangle\ge 1/n$. 
\end{lemma}

\proof
We fix $v\in S^{n-1}$ and define $S_+:=\{u\in S^{n-1}:\langle u,v\rangle\ge 0\}$ and $ S_-:=S^{n-1}\setminus S_+$. 
Since $\mu$ is centred and $\langle u,v\rangle\ge -1$, we have 
$$
-\int_{S_+}\langle u,v\rangle\, d\mu(u)=\int_{S_-}\langle u,v\rangle\, d\mu(u)\ge -\mu(S_-),
$$
and hence
\begin{equation}\label{mus-}
\mu(S_-)\ge \int_{S_+}\langle u,v\rangle\, d\mu(u).
\end{equation}
Choose $u^*\in \text{supp}\ \mu$ such that $\langle u^*,v\rangle=\max\{\langle u,v\rangle:u\in \text{supp}\ \mu\}$. The maximum exists as $\text{supp}\ \mu$ is compact. 
It is also clear (since $\mu$ is centred) that $u^*\in S_+$. Then \eqref{mus-} implies
\begin{align}\label{upperbound2}
\int_{S_+}\langle u,v\rangle^2\, d\mu(u)&\le \langle u^*,v\rangle\int_{S_+}\langle u,v\rangle\, d\mu(u)\le \langle u^*,v\rangle \mu(S_-).
\end{align}
In addition, we have
\begin{align}\label{upperbound3}
\int_{S_-}\langle u,v\rangle^2\, d\mu(u)&\le\int_{S_-}|\langle u,v\rangle|\, d\mu(u)=-\int_{S_-}\langle u,v\rangle\, d\mu(u)
=\int_{S_+}\langle u,v\rangle\, d\mu(u)\nonumber\\
&\le \langle u^*,v\rangle \mu(S_+).
\end{align}
Using \eqref{upperbound2},  \eqref{upperbound3}, the isotropy of $\mu$ and $\mu(S^{n-1})=n$, we conclude 
\begin{align*}
1&=\int_{S^{n-1}}\langle u,v\rangle^2\, d\mu(u)=\int_{S_+}\langle u,v\rangle^2\, d\mu(u)+\int_{S_-}\langle u,v\rangle^2\, d\mu(u)\\
&\le \langle u^*,v\rangle \mu(S_-)+\langle u^*,v\rangle \mu(S_+)=\langle u^*,v\rangle \mu(S^{n-1})=n\langle u^*,v\rangle,
\end{align*}
which yields the assertion.
\proofbox

We say that a non-empty closed subset $X$ of $S^1$ is proper, if for any $v\in S^1$, there exists some $u\in X$ 
such that $\langle v,u\rangle \geq \frac12$. A closed set $X\subset S^1$ is proper if and only if 
the angle of two consecutive points of $X$ is at most $2\pi/3$.

For a non-empty closed set $X\subset S^1$, let $d_0(X)$ be the minimum of $\delta_H(X,\sigma)$ where $\sigma$ runs through the set of contact points of the regular triangles circumscribed about $B^2$. If $X$ is proper, then clearly $d_0(X)\leq \pi/3$.

\begin{lemma}
\label{planarapprox}
If $X\subset S^1$ is proper, and $d_0(X)\geq\eta$ for $\eta\in(0,\frac{\pi}6]$, then there exist  $u,v\in X$ such that
$\eta\leq \angle(u,v)\leq \frac{2\pi}3-\eta$.
\end{lemma}
\proof  We prove the lemma by contradiction, thus we suppose that for any
$u,v\in X$, we have
\begin{equation}
\label{indirect-planar}
\mbox{either $\angle(u,v)<\eta$ or $\angle(u,v)>\frac{2\pi}3-\eta\geq  \frac{\pi}2>2\eta$}.
\end{equation}
The set $X$ has at least four elements since $X$ is proper and $d_0(X)>0$. Thus there exist $u'_ 1,v'_1\in X$ such that
$0<\angle(u'_1,v'_1)\leq  \frac{\pi}2$. We deduce from
(\ref{indirect-planar}) that
$\angle(u'_1,v'_1)<\eta$.  According to (\ref{indirect-planar}), there exists $v_1\in X$ such that $\angle(u'_1,v_1)$ is maximal under the conditions $\angle(u'_1,v_1)<\eta$ and $v'_1\in{\rm pos}\{u'_1,v_1\}$. Similarly, there exists $u_1\in X$ such that $\angle(u_1,v_1)$ is maximal under the conditions $\angle(u_1,v_1)<\eta$ and $u'_1\in{\rm pos}\{u_1,v_1\}$.

As $X$ is proper, there exists $u_2\in X$ such that ${\rm lin}\,v_1$ separates $u_1$ and $u_2$, and
$\angle(u_2,v_1)$ is minimal under the conditions $\angle(u_2,v_1)\leq \frac{2\pi}3$ and that 
${\rm lin}\,v_1$ separates $u_1$ and $u_2$. We actually have
\begin{equation}
\label{planaru2v1}
\mbox{$\frac{\pi}2\leq \frac{2\pi}3-\eta<\angle(u_2,v_1)\leq \frac{2\pi}3,$}
\end{equation}
since $\angle(u_2,v_1)<\eta$ would imply $\eta\leq \angle(u_2,u_1)<2\eta$, contradicting
(\ref{indirect-planar}). In particular, we have $X\cap \text{pos}\{u_2,v_1\}=\{u_2,v_1\}$. 
Similarly, there exists $v_3\in X$ such that  ${\rm lin}\,u_1$ separates $v_1$ and $v_3$, and
\begin{equation}
\label{planarv3u1}
\mbox{$\frac{\pi}2\leq\frac{2\pi}3-\eta<\angle(v_3,u_1)\leq \frac{2\pi}3,$}
\end{equation}
moreover $X\cap{\rm pos}\{v_3,u_1\}=\{v_3,u_1\}$. It also follows from (\ref{planaru2v1}) and 
(\ref{planarv3u1}) that $u_2$ and $v_3$ are not opposite, and the shorter arc of $S^1$ connecting 
them does not contain $u_1$ and $v_1$.

Finally, let $v_2\in X\cap {\rm pos}\{u_2,v_3\}$ maximize $\angle(v_2,u_2)$ under the condition
$\angle(v_2,u_2)<\eta$, and let $u_3\in X\cap {\rm pos}\{u_2,v_3\}$ maximize $\angle(u_3,v_3)$ under the condition
$\angle(u_3,v_3)<\eta$. Here possibly $v_2=u_2$ or $u_3=v_3$. If there were
$w\in  X\cap{\rm int}\,{\rm pos}\{v_2,u_3\}$, then $\angle(w,v_3)>\frac{\pi}2$ and $\angle(w,u_2)>\frac{\pi}2$
would follow from (\ref{indirect-planar}), what is absurd. Therefore $X\cap{\rm pos}\{u_3,v_2\}=\{u_3,v_2\}$, and
\begin{equation}
\label{planaru2v3}
\mbox{$\frac{\pi}2\leq\frac{2\pi}3-\eta<\angle(u_3,v_2)\leq \frac{2\pi}3,$}
\end{equation}

Now the arcs $S^1\cap {\rm pos}\{u_1,v_2\}$, $S^1\cap {\rm pos}\{u_2,v_3\}$ and $S^1\cap {\rm pos}\{u_3,v_1\}$
cover $S^1$ by their constructions, thus
\begin{equation}
\label{planararcsum}
\angle(u_1,v_2)+\angle(u_2,v_3)+\angle(u_3,v_1)>2\pi.
\end{equation}
In particular, one of $\angle(u_1,v_2)$, $\angle(u_2,v_3)$ and $\angle(u_3,v_1)$  is larger than $\frac{2\pi}3$ by
(\ref{planararcsum}).

If $\angle(u_1,v_2)>\frac{2\pi}3$, then we define $p_3\in S^1$ in such a way that $-p_3$ is the midpoint of the arc
$S^1\cap {\rm pos}\{u_1,v_2\}$. For $i=1,2$, let $p_i\in S^1$ satisfy $\angle(p_i,p_3)=\frac{2\pi}3$ in such a way 
that $p_1$ and $p_2$ lie on the same side of ${\rm lin}\ p_3$ where $u_1$ and $v_2$ lie, respectively. In particular,
$p_1$, $p_2$ and $p_3$ are vertices of a regular triangle. We deduce using (\ref{planarv3u1}) and
(\ref{planaru2v3})  that
\begin{equation}
\label{planarp1p2}
p_1,p_2\in {\rm pos}\{u_1,v_2\}\mbox{ \ and \
$\angle(u_1,v_2)<\frac{2\pi}3+2\eta$.}
\end{equation}
For $i=1,2$, it follows from (\ref{planarp1p2}) that if $w\in S^1\cap {\rm pos}\{u_i,v_i\}$, then
$\angle(w,p_i)<\eta$. In addition, (\ref{planarv3u1}) and
(\ref{planaru2v3}) yield that  if $w\in S^1\cap {\rm pos}\{u_3,v_3\}$, then
$\angle(w,p_1)<\eta$, and hence $d_0(X)<\eta$, which is a contradiction. If
$\angle(u_2,v_3)>\frac{2\pi}3$ or $\angle(u_3,v_1)>\frac{2\pi}3$ in
(\ref{planararcsum}), then similar arguments lead to a contradiction, which completes the proof of Lemma~\ref{planarapprox}.
\proofbox

In the following, we use the fact (T) that for $0\le \beta\le\alpha\le 2\pi/3$ the function
$$
F(t)=\tan\left(\frac{\alpha+t}{2}\right)+\tan\left(\frac{\beta-t}{2}\right)=\frac{2\sin\left(\frac{\alpha+\beta}{2}\right)}{\cos\left(\frac{\alpha+\beta}{2}\right)+\cos\left(t+\frac{\alpha-\beta}{2}\right)}
$$
is increasing for $0\le t\le \min\{\beta,\frac{2\pi}{3}-\alpha\}$. 

\bigskip

After these preparations, we turn to the proof of Theorem \ref{planariso}. 

\proof It is sufficient to  prove that if  $\eta\in(0,\frac{\pi}6]$, and
$d_0({\rm supp}\,\mu)\geq \eta$, then
\begin{equation}
\label{planarvolbig}
V(Z(\mu))\leq \left(1-\frac{\eta}{8}\right)V(T^2).
\end{equation}
Indeed, if $d_0({\rm supp}\,\mu) >32\varepsilon$, then $8\varepsilon<\pi/6$, since $d_0({\rm supp}\,\mu)\le 2\pi/3$ by Lemma \ref{lowerbound}. But then the 
preceding claim can be applied with $\eta=8\varepsilon$.

Now we turn to the proof of the claim. 
It follows from Lemma~\ref{planarapprox} that there exist $u_1,u_2\in {\rm supp}\,\mu$ such that
\begin{equation}
\label{planarvolangle}
\mbox{$\eta\leq \angle(u_1,u_2)\leq \frac{2\pi}3-\eta$.}
\end{equation}

 Since by Lemma \ref{lowerbound} ${\rm supp}\,\mu$ is proper, there exist 
$u_3,\ldots,u_k\in {\rm supp}\,\mu$, $k\geq 4$, such that $u_1,\ldots,u_k$  (in this order) lie on $S^1$ and form a proper set. 
Then 
$$
V(Z(\mu))\le 2\sum_{i=1}^k\tan\left(\frac{\alpha_i}{2}\right),
$$
where $\alpha_1=\angle(u_1,u_2)\in [\eta,\frac{2\pi}{3}-\eta]$, $\alpha_i=\angle(u_i,u_{i+1})$ with $u_{k+1}:=u_1$ and $0\le\alpha_i\le 2\pi/3$. 
Applying repeatedly (T) to pairs of the angles $\alpha_2,\ldots,\alpha_k$, it follows that 
\begin{align*}
2\sum_{i=1}^k\tan\left(\frac{\alpha_i}{2}\right)&=2\left(\tan\left(\frac{\alpha_1}{2}\right)+
\tan\left(\frac{2\pi-\frac{4\pi}{3}-\alpha_1}{2}\right)+2\tan\left(\frac{\pi}{3}\right)
\right)\\
&=2\left(\tan\left(\frac{\alpha_1}{2}\right)+\tan\left(\frac{\pi}{3}-\frac{\alpha_1}{2}\right)+2\sqrt{3}
\right)\\
&\le 2\left(\tan\left(\frac{\eta}{2}\right)+\tan\left(\frac{\pi}{3}-\frac{\eta}{2}\right)+2\sqrt{3}
\right)\\
&\le 2\left(\frac{\sqrt{3}}{\frac{1}{2}\cos\left(\eta-\frac{\pi}{3}\right)}+2\sqrt{3}\right)\le \left(1-\frac{\eta}{8}\right)6\sqrt{3}\\
&=\left(1-\frac{\eta}{8}\right)V(T^2),
\end{align*}
which proves the assertion.
\proofbox

\section{Proof of Theorem \ref{voln=2}}\label{sec12}
Let $K$ be a convex body in $\R^2$ whose John ellipsoid is the Euclidean unit ball. As before (at the beginning of Section \ref{secisoperimetric}), the contact points of $K$ and $B^2$ define a discrete, centred, isotropic measure $\mu$ and a polytope $Z=Z(\mu)$ which contains $K$. 

If $V(Z)\ge (1-\varepsilon)V(T^2)$ with some $\varepsilon\in (0,1)$, then Theorem \ref{planariso} implies the existence of a regular simplex $S$ circumscribed about $B^2$ such that $\delta_{H}({\rm supp}\ \mu,{\rm supp}\ \mu_S)\le 32\ \varepsilon$. Choosing $\eta:=32\ \varepsilon<1/18$, that is with $\varepsilon<1/(18\cdot 32)$, we see from 
Lemma \ref{closesimplex} that $d(Z)<18\cdot 32\ \varepsilon$. Hence, if $d(Z)\ge 18\cdot 32\ \varepsilon$ and $\varepsilon<1/(18\cdot 32)$, then $V(Z)< (1-\varepsilon)V(T^2)$, and therefore $S(K)^2/V(K)\le (1-\varepsilon){\rm ir}(T^2)$. On the other hand, if $d(Z)<18\cdot 32\ \varepsilon$ and $\delta_{\rm vol}(K,T^2)\ge 16\cdot 18\cdot 32\  \varepsilon$, then Lemma \ref{ZK} (ii) 
implies that 
$$
\frac{S(K)^2}{V(K)}\le \left(1-\frac{1}{8}\ 16\cdot 18\cdot 32\ \varepsilon\right){\rm ir}(T^2)=(1- 16\cdot 32\ \varepsilon){\rm ir}(T^2),
$$
provided that $16\cdot 18\cdot 32 \ \varepsilon<1$. This implies the assertion of the theorem.
\proofbox

\bigskip

\noindent{\bf Acknowledgement: }  We are grateful to Rolf Schneider for initiating the problem, and  to Erwin Lutwak, Vitali Milman and Gaoyong Zhang for providing insight into properties of isotropic measures. Special thanks are due to Keith Ball for sharing a concise proof of Lemma~\ref{ciuibig}.

\bigskip

Authors' addresses:

\bigskip

K\'aroly J. B\"or\"oczky,
Alfr\'ed R\'enyi Institute of Mathematics, Hungarian Academy
of Sciences 1053 Budapest, Re\'altanoda u.~13-15. HUNGARY.
E-mail: carlos@renyi.hu
\medskip

Daniel Hug, 
Karlsruhe Institute of Technology (KIT),
D-76128 Karlsruhe, Germany. E-mail: daniel.hug@kit.edu

\end{document}